\documentclass{gtart_a}
\pdfoutput=1

\usepackage[matrix,arrow,curve]{xy}
\usepackage{amscd,graphicx}

\title[Nonstabilized Nielsen coincidence invariants]{Nonstabilized
Nielsen coincidence invariants and Hopf--Ganea homomorphisms}

\author{Ulrich Koschorke}
\givenname{Ulrich}
\surname{Koschorke}
\address{Universit\"at Siegen\\
Emmy Noether Campus\\
Walter-Flex-Str. 3\\\newline
D-57068 Siegen\\
Germany}
\email{koschorke@mathematik.uni-siegen.de}
\urladdr{http://www.math.uni-siegen.de/topology/}

\volumenumber{10}
\issuenumber{}
\publicationyear{2006}
\papernumber{17}
\lognumber{0662}
\startpage{619}
\endpage{666}

\doi{}
\MR{}
\Zbl{}


\received{29 September 2005}
\revised{9 March 2006}
\accepted{21 April 2006}
\proposed{Shigeyuki Morita}
\seconded{Wolfgang L\"uck, Tom Goodwillie}
\published{24 May 2006}
\publishedonline{24 May 2006}
\corresponding{}
\editor{}
\version{}

\let\xysavmatrix\xymatrix
\def\xymatrix{\disablesubscriptcorrection\xysavmatrix}
\AtBeginDocument{\let\bar\wbar\let\tilde\wtilde}

\numberwithin{equation}{section}
\newcommand{\we}{\smash{\rlap{\kern 6pt 
\raise 4pt\hbox{\footnotesize $\sim$}}}\longrightarrow}

\makeatletter
\def\cnewtheorem#1[#2]#3{\newtheorem{#1}{#3}[section]
\expandafter\let\csname c@#1\endcsname\c@thm}

\newtheorem{thm}{Theorem}[section]
\cnewtheorem{prop}[thm]{Proposition}
\cnewtheorem{cor}[thm]{Corollary}
\cnewtheorem{lem}[thm]{Lemma}
\newtheorem*{lemA}{Lemma A}

\theoremstyle{definition}
\cnewtheorem{definition}[thm]{Definition}
\cnewtheorem{example}[thm]{Example}
\newtheorem*{question}{Question}
\cnewtheorem{notcon}[thm]{Conventions}
\newtheorem*{stepI}{Step I}
\newtheorem*{stepII}{Step II}

\theoremstyle{remark}
\cnewtheorem{remark}[thm]{Remark}

\makeatother  

\newcommand{\att}{\operatorname{att}}
\newcommand{\pinch}{\operatorname{pinch}}
\newcommand{\proj}{\operatorname{proj}}
\newcommand{\pr}{\operatorname{pr}}
\newcommand{\incl}{\operatorname{incl}}
\newcommand{\quot}{\operatorname{quot}}
\newcommand{\coll}{\operatorname{coll}}
\newcommand{\rel}{\operatorname{rel}}

\newcommand{\dist}{\operatorname{dist}}

\newcommand{\id}{\operatorname{id}}
\newcommand{\inv}{\operatorname{inv}}
\newcommand{\im}{\operatorname{im}}
\newcommand{\In}{\operatorname{in}}
\newcommand{\ind}{\operatorname{ind}}
\newcommand{\forg}{\operatorname{forg}}
\newcommand{\const}{\operatorname{const}}

\newcommand{\stab}{\operatorname{stab}}
\newcommand{\ret}{\operatorname{ret}}
\newcommand{\scirc}{{\scriptstyle{\circ}}}

\newcommand{\MCC}{\mathit{MCC}}
\newcommand{\MC}{\mathit{MC}}
\newcommand{\MN}{\mathit{MN}}

\begin{document}

\begin{asciiabstract}
In classical fixed point and coincidence theory the notion of Nielsen
numbers has proved to be extremely fruitful. We extend it to pairs
(f_1,f_2) of maps between manifolds of arbitrary dimensions, using
nonstabilized normal bordism theory as our main tool. This leads to
estimates of the minimum numbers MCC(f_1,f_2) (and
MC(f_1,f_2), respectively) of path components (and of points, resp.) in
the coincidence sets of those pairs of maps which are homotopic to
(f_1,f_2). Furthermore, we deduce finiteness conditions for
MC(f_1,f_2). As an application we compute both minimum numbers
explicitly in various concrete geometric sample situations.

The Nielsen decomposition of a coincidence set is induced by the
decomposition of a certain path space E(f_1,f_2) into
path components. Its higher dimensional topology captures further
crucial geometric coincidence data. In the setting of homotopy groups
the resulting invariants are closely related to certain Hopf--Ganea
homomorphisms which turn out to yield finiteness obstructions for
MC.
\end{asciiabstract}

\begin{htmlabstract}
<p class="noindent">
In classical fixed point and coincidence theory the notion of
Nielsen numbers has proved to be extremely fruitful. We extend
it to pairs (f<sub>1</sub>,f<sub>2</sub>) of maps between
manifolds of arbitrary dimensions, using nonstabilized normal
bordism theory as our main tool. This leads to estimates of the
minimum numbers <em>MCC</em>(f<sub>1</sub>,f<sub>2</sub>) (and
<em>MC</em>(f<sub>1</sub>,f<sub>2</sub>), resp.) of path components (and
of points, resp.) in the coincidence sets of those pairs of maps which
are homotopic to (f<sub>1</sub>,f<sub>2</sub>). Furthermore, we deduce
finiteness conditions for <em>MC</em>(f<sub>1</sub>,f<sub>2</sub>). As
an application we compute both minimum numbers explicitly in various
concrete geometric sample situations.</p>

<p class="noindent">
The Nielsen decomposition of a coincidence set is induced by the
decomposition of a certain path space E(f<sub>1</sub>,f<sub>2</sub>)
into path components. Its higher dimensional topology captures further
crucial geometric coincidence data. In the setting of homotopy groups
the resulting invariants are closely related to certain Hopf&ndash;Ganea
homomorphisms which turn out to yield finiteness obstructions for
<em>MC</em>.</p>
\end{htmlabstract}

\begin{abstract}
In classical fixed point and coincidence theory the notion of Nielsen
numbers has proved to be extremely fruitful. We extend it to pairs
$(f_1,f_2)$ of maps between manifolds of arbitrary dimensions, using
nonstabilized normal bordism theory as our main tool. This leads to
estimates of the minimum numbers $\mathit{MCC}(f_1,f_2)$ (and
$\mathit{MC}(f_1,f_2)$, resp.) of path components (and of points, resp.) in
the coincidence sets of those pairs of maps which are homotopic to
$(f_1,f_2)$. Furthermore, we deduce finiteness conditions for
$\mathit{MC}(f_1,f_2)$. As an application we compute both minimum numbers
explicitly in various concrete geometric sample situations.

The Nielsen decomposition of a coincidence set is induced by the
decomposition of a certain path space $E(f_1,f_2)$ into
path components. Its higher dimensional topology captures further
crucial geometric coincidence data. In the setting of homotopy groups
the resulting invariants are closely related to certain Hopf--Ganea
homomorphisms which turn out to yield finiteness obstructions for
$\mathit{MC}$.
\end{abstract}

\maketitle

\section{Introduction}
\label{sec1}

In this paper we develop a coherent geometric approach to coincidence
phenomena in arbitrary codimensions. We prove (and extend considerably)
results which were announced in part in \cite{K5} and \cite{K6}.

Consider two continuous maps $f_1, f_2 \co  M \to N$ between smooth connected
manifolds without boundary, of arbitrary positive  dimensions $m$ and $n,
M$ being compact.

We would like to measure how small (or simple in some sense) the
coincidence locus
\begin{equation}   
\label{1.1}
C (f_1, f_2) := \{ x \in M | f_1 (x) = f_2 (x) \}
\end{equation}
can be made by varying $f_1$ and $f_2$ within their homotopy classes.

One possible measure is the classical {\bf m}inimum number of {\bf
c}oincidence points
\begin{equation}              
\label{1.2}
\MC (f_1, f_2) := \min \{ \# C (f'_1, f'_2) | f'_1 \sim f_1, f'_2
\sim f_2\}
\end{equation}
(cf Bogaty\u\i--Gon\c{c}alves--Zieschang \cite[1.1]{BGZ}). It coincides with
the minimum number $\min \{ \# C (f'_1, f_2) | f'_1 \sim f_1\}$ where only
$f_1$ is modified by a homotopy (cf Brooks \cite{Br}). In particular, in
topological fixed point theory (where $M = N$ and $f_2 =$ identity) this
minimum number is the principal object of study (cf Brown \cite[page~9]{B}).

However, in higher codimensions the coincidence locus is generically
a manifold of dimension $m-n > 0$, and $\MC (f_1, f_2)$ is often
infinite (see eg Examples \ref{1.12}, \ref{1.25}, \ref{5.6}, and
\ref{6.21} below). Thus
in many situations it seems more meaningful to study the {\bf m}inimum
number of {\bf c}oincidence {\bf c}omponents
\begin{equation}                                    
\label{1.3}
\begin{split}
\MCC (f_1, f_2)
 :=& \min \{ \# \pi_0 (C (f'_1, f'_2)) | f'_1 \sim f_1, f'_2 \sim
 f_2\} \\
 =& \min \{ \# \pi_0 (C (f'_1, f_2)) | f'_1 \sim f_1 \}
\end{split}
\end{equation}
where $\# \pi_0 (C (f'_1, f'_2))$ denotes the (generically finite) number
of {\em path components} of the indicated coincidence subspace of $M$
(compare Bogaty\u\i--Gon\c{c}alves--Zieschang \cite[page~47, line~3]{BGZ}).

\begin{definition}                           
\label{1.4}
The pair of maps $(f_1, f_2)$ is called {\it loose} if $\MC (f_1, f_2)
= 0$
(or, equivalently, $\MCC (f_1, f_2) = 0)$, ie if the maps $f_1$
and $f_2$ can be deformed away from one another.
\end{definition}

\begin{question}
How big are $\MCC (f_1, f_2)$ and $\MC (f_1, f_2)$? In particular, when
do these invariants vanish, ie when is the pair $(f_1, f_2)$ loose?
\end{question}

In order to attack this problem let us study the geometry of generic
coincidence submanifolds.

After performing an approximation we may assume that the map $(f_1, f_2)
\co  M \to N \times N$ is smooth and transverse to the diagonal $\Delta = \{
(y, y) \in N \times N | y \in N\}$.

Then the coincidence locus
\begin{equation}                           
\label{1.5}
C = C (f_1, f_2) = (f_1, f_2)^{-1} (\Delta) = \{ x \in M | f_1  (x) =
f_2 (x) \}
\end{equation}
is a closed smooth $(m{-}n)$--dimensional submanifold of $M$. It comes
with two important data. First there is a commuting diagram of maps
\begin{equation}                                         
\label{1.6}
  \xymatrix{
  & E(f_1,f_2) \ar[d]^-{\text{pr}}&
  \!\!\!\!\!\!\!\!\!\!\!\!\!\!\!\!
  :=
  \{(x,\theta)\in M\times P(N) \mid \theta(0)=f_1(x); \theta(1)=f_2(x)\}
  \\
  C \ar[ur]^-{\wtilde{g}}
  \ar[r]_-{g=\text{incl}} & M
  }
\end{equation}\eject
where $P(N)$ (and $\pr$, resp.), denote the space of all continuous
paths $\theta \co  [0, 1] \to N$, endowed with the compact--open topology,
(and the obvious projection, resp.);  the lifting $\wtilde g$ adds
the constant path at $f_1 (x) = f_2 (x)$ to $g (x) = x \in C$. The second
datum is the (composite) vector bundle isomorphism
\begin{equation}                                       
\label{1.7}
\wbar g^\# \co  \nu (C, M) \cong ((f_1, f_2) | C)^* (\nu (\Delta,
N \times N)) \cong f_1^* (TN) | C
\end{equation}
which describes the normal bundle of $C$ in $M$ (see the figure in
Koschorke \cite[Section~4]{K6} for an illustration).

The resulting bordism class
\begin{equation}                                          
\label{1.8}
\omega^\# (f_1, f_2) = [C (f_1, f_2), \wtilde g, \wbar g^\# ] \
\in \Omega^\# (f_1, f_2)
\end{equation}
in an appropriate bordism set (cf \eqref{2.1} and \eqref{2.2} below) is our key
coincidence invariant. It turns out that the lifting $\wtilde g$
plays a crucial role. Indeed, in general the path space $E (f_1, f_2)$
(cf \eqref{1.6}) has a very rich topology involving both $M$ and the loop
space of $N$ (cf Koschorke \cite[2.1]{K3}). Already the set $\pi_0 (E (f_1,
f_2))$ of path components can be huge -- it corresponds bijectively to
the Reidemeister set
\begin{equation}                               
\label{1.9}
R (f_1, f_2) = \pi_1 (N) / \text{Reidemeister equivalence}
\end{equation}
(compare Bogaty\u\i--Gon\c{c}alves--Zieschang \cite[3.1]{BGZ} and Koschorke
\cite[2.1]{K3}) which is of central
importance in classical Nielsen theory. This leads to a natural
decomposition
\begin{equation*}
C (f_1, f_2) =  \coprod_{A \in \pi_0 (E (f_1, f_2))} \wtilde g^{-
1} (A) .
\end{equation*}
We define $N^\# (f_1, f_2)$ to be the corresponding number of nontrivial
contributions by the various path components $A$ of $E (f_1, f_2)$ to
$\omega^\# (f_1, f_2)$ (see \fullref{2.7} below).

If we forget the fact that the manifold $C (f_1, f_2)$ is {\it embedded}
in $M$ and if we stabilize $\wbar g^\#$ to yield only a description
of the {\it stable} normal bundle of $C (f_1, f_2)$ we obtain the normal
bordism class
\begin{equation}                              
\label{1.10}
\wtilde\omega (f_1, f_2) = [C (f_1, f_2), \wtilde g, \wbar g]
\in \Omega_{m-n} (E (f_1, f_2); \wtilde\varphi := \pr^*
(f^*_1 (TN) - TM))
\end{equation}
and the corresponding Nielsen number $N (f_1, f_2)$. These \lq\lq
stabilized\rq\rq\ invariants were studied in detail in \cite{K3}.

Let us put our approach into perspective. Recall the decisive progress
made by J Nielsen on the classical minimizing problem when he decomposed
fixed point sets into equivalence classes.  In our interpretation this is
just the decomposition of a 0--dimensional bordism class according to the
path components of its target space. In higher (co)dimensions $(m{-}n)$ the
map $\wtilde g$ into $E (f_1, f_2)$ and the \lq\lq twisted
framing\rq\rq\  $\wbar g^\#$ contain much more information and lead
sometimes to a complete calculation of $\MCC (f_1, f_2)$ and $\MC (f_1,
f_2)$ (cf, for example, \fullref{1.25} below or Koschorke \cite[Examples~I--IV]{K6}).

\begin{thm}                                        
\label{1.11}
{\rm (i)}\qua The Nielsen numbers $N (f_1, f_2)$ and $N^\# (f_1, f_2)$
are finite and depend only on the homotopy classes of $f_1$ and $f_2$;

{\rm (ii)}\qua $N (f_1, f_2) = N (f_2, f_1)$ and $N^\# (f_1, f_2) = N^\#
(f_2, f_1)$;

{\rm (iii)}\qua $0 \le N (f_1, f_2) \le N^\# (f_1, f_2) \le \MCC (f_1, f_2)
\le \MC (f_1, f_2)$; if $n \ne 2$, then also
$$\MCC (f_1, f_2) \le \# \pi_0 (E (f_1, f_2));$$
if $(m, n) \ne (2,2)$, then
$$\MC (f_1, f_2) \le \# \pi_0 (E (f_1, f_2))\quad\text{or}\quad
\MC (f_1, f_2) = \infty;$$
{\rm (iv)}\qua if $m = n$, then $N (f_1, f_2) = N^\# (f_1, f_2)$ coincides
with the classical Nielsen number (cf Bogaty\u\i--Gon\c{c}alves--Zieschang
\cite[Definition~3.6]{BGZ}).
\end{thm}

The proof and further details concerning our $\omega$--invariants and
Nielsen numbers will be given in \fullref{sec2} below. \fullref{sec3} is dedicated to
the minimum number $\MC (f_1, f_2)$.

\begin{example}                                  
\label{1.12}
Assume $N = S^1$. Then both Nielsen numbers of $(f_1, f_2)$ agree with
$\MCC (f_1, f_2)$ and are characterized by the identity
\begin{equation*}
(f_{1*} - f_{2*}) (H_1 (M; \mathbb Z)) = N^{(\#)} (f_1, f_2) \cdot H_1
(S^1; \mathbb Z) .
\end{equation*}
If $f_1 \sim f_2$, then $\# \pi_0 (E (f_1, f_2)) = \infty$ and
\begin{equation*}
N (f_1, f_2) = N^\# (f_1, f_2) = \MCC (f_1, f_2) = \MC (f_1, f_2) = 0 .
\end{equation*}
If $f_1$ and $f_2$ are not homotopic then
\begin{equation*}
N (f_1, f_2) = N^\# (f_1, f_2) = \MCC (f_1, f_2) = \# \pi_0 (E (f_1,
f_2)) \ne 0
\end{equation*}
and
\begin{equation*}
\MC (f_1, f_2) = \begin{cases}
N (f_1, f_2) & \text{if } m = 1; \\
\infty & \text{if } m \ge 2
\end{cases}
\end{equation*}
(Clearly in all other cases where $m = 1$ or $n = 1$ we have
$$N (f_1, f_2) = N^\# (f_1, f_2) = \MCC (f_1, f_2) = \MC (f_1, f_2) = 0.)$$
\end{example}

In higher codimensional coincidence theory two settings are of particular
interest. In one of them, the so-called root case (cf
Bogaty\u\i--Gon\c{c}alves--Zieschang \cite[page~69]{BGZ}) $f_2 = *$
is constant (henceforth our notation will not distinguish
between constant maps and their values). Here our invariants yield the
\lq\lq degrees\rq\rq\
\begin{equation}                                        
\label{1.13}
\deg^\# (f) = \omega^\# (f, *), \quad \widetilde\deg (f) =
\wtilde\omega (f, *)
\end{equation}
of a given map $f \co  M \to N$. The choice of the constant $* \in N$ is
not truly significant since any path joining two such constants induces
a bijection between the corresponding bordism sets which is compatible
with degrees and Reidemeister decompositions. Note also that $E (f, *)$
is the mapping fiber of $f$ (cf Whitehead \cite[I.7]{W})

The second particularly interesting setting concerns selfcoincidences
(where $f_1$ is equal or at least homotopic to $f_2$). Here we know from
the very outset that $\MCC (f_1, f_2) \le 1$ (since $C (f_1, f_1) = M)$.
The remaining question whether $\MCC (f_1, f_2) = 0$ or, equivalently,
whether $f_1$ can be deformed away from itself was studied in \cite{K2}
(and related -- in one particular example -- to a fascinating problem
concerning Lie groups and their role in homotopy theory). It is also
worthwhile noting that each of the selfcoincidence invariants $\omega^\#
(f, f)$ and $\wtilde\omega (f, f)$ is determined by the corresponding
degree (cf \eqref{5.3} and \fullref{5.5} below).

Now let a map $f \co  M \to N$ and a constant $* \in N$ be given.

\begin{thm}                        
\label{1.14}
Consider the pairs $(f_1, f_2) = (f, *)$ (root case) and $(f_1, f_2)
= (f, f)$ (selfcoincidence case) simultaneously. Define $b (f_1,
f_2)$ by $b (f, *) := \# \pi_0 (E (f, *))$  = index of $f_* (\pi_1 (M))$
in $\pi_1 (N)$,  and $b (f, f) := 1$.

In both cases the following holds:

{\rm (i)}\qua If $\omega^\# (f_1, f_2) \ne 0$ then $N^\# (f_1, f_2) = b
(f_1, f_2)$; if in addition $n \ne 2$ then this Nielsen number agrees with
$\MCC (f_1, f_2)$ (and also with the minimum number $\MC (f_1, f_2)$
whenever it is finite).

{\rm (ii)}\qua If even $\wtilde\omega (f_1, f_2) \ne 0$ then also $N
(f_1, f_2) = N^\# (f_1, f_2)$.
\end{thm}

In particular in both cases the Nielsen numbers $N^\# (f_1, f_2)$ and $N
(f_1, f_2)$ can take only the values $0$ and $b (f_1, f_2)$.

This has important consequences. First of all since $b (f, *)$ depends
on $f$ but $b (f, f)$ does not, the compatibilities with covering spaces
must be different in the root and selfcoincidence settings.

\begin{prop}                          
\label{1.15}
If $\wtilde f \co  M \to \wtilde N$ is a lifting of $f$ to any
connected $d$--fold covering space $\wtilde N$ of $N$ and $\wtilde*
\in \wtilde N, * \in N$ are arbitrary constants, then
$$N^\# (f, *) = d \cdot N^\# (\wtilde f, \wtilde*) \quad
\text{but } N^\# (f, f) = N^\# (\wtilde f, \wtilde f)$$
and the corresponding identities hold also for the weaker (stabilized)
Nielsen numbers $N (f, *)$ and $N (f, f)$.
\end{prop}

(This holds even when $d = \infty$, provided we define $\infty \cdot 0
= 0$).

\begin{remark}                            
\label{1.16}
If $\deg^\# (f)$ and $\omega^\# (f, f)$, resp., do not vanish and $n \ne
2$, then these identities still hold when we replace $N^\#$ by $\MCC$
(and by $\MC$, provided $\MC (f, *)$ and $\MC (f, f)$, resp., are finite;
compare \fullref{1.14}). In general, however, we can only establish the inequalities
\begin{equation*}
\MC (f, *) \ge d \cdot \MC (\wtilde f, \wtilde *) \quad
\text{and} \quad \MC (f, f) \ge \MC (\wtilde f, \wtilde f) .
\end{equation*}
At least in the root case the corresponding {\em equality} is often not
valid. Indeed there are many examples where $\deg^\# (f) \ne 0$ and $\MC
(\wtilde f, \wtilde *) = 1 $ while $\MC (f, *)$ (but not $d$)
is infinite (see Examples~\ref{1.25} and~\ref{6.21} below). Such phenomena are closely
related to Wyler's theory of injective points (cf Wyler \cite{Wy} and our
discussion following \fullref{6.18}).
\end{remark}

\fullref{1.14}, \fullref{1.15}, and further results concerning the root
and selfcoincidence settings will be proved in Sections~\ref{sec4}
and~\ref{sec5} below
(cf \fullref{4.17} and the discussions following \fullref{4.18} and
\eqref{5.3}).

As an illustration we test our approach in \fullref{sec6} in the case where $M$
is a sphere. Here
we can exploit two important advantages. On one hand there is a natural
identification of the bordism set $\Omega^\# (f_1, f_2)$ with a fixed
{\em group} which does not vary with $f_1$ and $f_2$. On the other hand
the algebraic structure of homotopy groups yields a certain homogeneity;
this allows considerable extensions of results which originally are
characteristic for the root setting. First we need to recall the
following notion.

\begin{definition}[cf Brown--Schirmer \cite{BS}]    
\label{1.17}
A map $f \co  M \to N$ is {\em not coincidence producing} if there exists
another map $\wbar f \co  M \to N$ such that the pair $(f, \wbar f)$
is loose (cf \fullref{1.4}).
\end{definition}

Such is always the case when $N$ allows a fixed point free selfmap $a\co  
N \to N$, eg when the manifold $N$ is open or its Euler number $\chi
(N)$ vanishes or $N = S^n$ or $N$ is the total space of a nontrivial
covering. (On the other hand one can easily exhibit settings where a vast
majority of maps is coincidence producing; see eg \eqref{6.19} below).

In the case $M = S^m$ non coincidence producing maps determine the
subgroup
\begin{equation}                                         
\label{1.18}
\pi^{(2)}_m (N) = \im p_* = \ker \partial \subset \pi_m (N)
\end{equation}
which arises also naturally in the exact homotopy sequence
\begin{equation}                                          
\label{1.19}
\begin{CD}
\cdots \longrightarrow \pi_m (\wtilde C_2 (N)) @>{ p_* }>> \pi_m (N)
@>{ \partial }>> \pi_{m -1} (N - \{ *\}) \longrightarrow \cdots
\end{CD}
\end{equation}
of the fibration $p  \co   \wtilde C_2 (N)  \to  N$ of the configuration
space of ordered pairs of distinct points in $N$.

\begin{thm}                         
\label{1.20}
Given $m \ge 1$, consider maps $f_i \co  S^m \to N$, $i = 1,2$, which
are not both coincidence producing. Assume $\omega^\# (f_1, f_2) \ne 0$.

Then
\begin{equation*}
N^\# (f_1, f_2) = \# \pi_0 (E (f_1, f_2)) = [\pi_1 (N): (f_{1 *} -
f_{2 *}) (\pi_1 (S^m))] .
\end{equation*}

If in addition $n \ne 2$, then this Nielsen number agrees with $\MCC (f_1,
f_2)$ (and also with the minimum number $\MC (f_1, f_2)$ whenever it
is finite).

If $\wtilde\omega (f_1, f_2) \ne 0$ and $n \ge 1$, then also $N (f_1,
f_2) = N^\# (f, f_2)$.
\end{thm}

\begin{cor}                             
\label{1.21}
If $f \co  S^m \to N$ is not coincidence producing and $\pi_1 (N) \ne 0$
then $\omega^\# (f, f) = 0$.
\end{cor}

Indeed otherwise $N^\# (f, f)$ would have to agree both with $\# \pi_1
(N)$ and with $1$. This highlights a rather astonishing feature of our
invariant in the non simply connected case: if $f$ can be deformed away
from  {\em any} map, $\omega^\#$ behaves as if $f$ can be deformed away
from itself.

\begin{cor}                                 
\label{1.22}
If $\pi_1 (N)$ has a nontrivial proper subgroup $G$ then $\omega^\#
(f, f) = 0$ for \emph{every} map $f \co  S^m \to N$.
\end{cor}

Indeed $G$ corresponds to a nontrivial covering space $\wtilde N$
of $N$ with $\pi_1 (\wtilde N) \ne 0$. Thus a lifting $\wtilde
f$ of $f$ is not coincidence producing and we have here $0 = N^\#
(\wtilde f, \wtilde f) = N^\# (f, f)$ (cf \fullref{1.15}).

For certain specific target manifolds $N$ (eg when $N$ is open or when
$N \ne S^1$ has an infinite fundamental group or when $N$ is a nontrivial
product of manifolds) it is rather easy to see that all pairs of maps
$f_1, f_2 \co  S^m \to N$ are loose (see \fullref{6.3} below). In view of
\fullref{1.20} a more systematic and detailed discussion is desirable: what
happens in general when $\omega^\# (f_1, f_2)$ vanishes? It is reasonable
to concentrate first on the root case. Given $m \ge 1$, we define (in
\fullref{6.4})
an abelian group $X_m (N)$ which measures to some extend the comparative
strength of the Nielsen number $N^\# (f, *)$ on one hand and of $\MCC (f,
*)$ on the other hand, taking into account all maps $f \co  S^m \to N$. This
group vanishes precisely if the conditions $N^\# (f, *) = 0$ and $\MCC
(f, *) = 0$ are equivalent. Actually this is often satisfied (cf
\fullref{6.6} below).

\begin{thm}                      
\label{1.23}
Given $m \ge 1$, assume $X_m (N) = 0$. Then a pair of maps $f_1, f_2\co  
S^m \to N$ is loose if and only if $\omega^\# (f_1, f_2) = 0$ and at
least one of the maps $f_1, f_2$ is not coincidence producing.
\end{thm}

Thus, if in the situation of \fullref{1.20} $X_m (N)$ and $\omega^\# (f_1,
f_2)$ vanishes, then so do also $N^\# (f_1, f_2), \MCC (f_1, f_2)$
and $\MC (f_1, f_2)$.

\begin{example}[$(N = S^n)$]                        
\label{1.24}
 Consider maps $f_1, f_2 \co  S^m \to S^n$ where $m, n \ge 1$, and let $a$
 denote the antipodal involution. Then
$$\MCC (f_1, f_2) = N^\# (f_1, f_2) = \
\begin{cases}
0 & \text{if } f_i \sim a f_2;\\
\#\pi_0(E(f_1,f_2)) & \text{otherwise.}
\end{cases}$$
If $f_1 \not\sim a f_2$ then $\# \pi_0 (E (f_1, f_2))$ equals $1$ (and
$| d^0 (f_1) - d^0 (f_2) |$, resp.) according as $n \ne 1$ (or $m = n =
1$, resp.; here $d^0 (f_i) \in \mathbb Z$ is the usual degree).

This shows that the \lq\lq strong\rq\rq\ Nielsen number $N^\# (f_1, f_2)$
(based on the nonstabilized invariant $\omega^\# (f_1, f_2))$ is often
strictly larger than $N (f_1, f_2)$.  For example, \cite[Corollary~1.17]{K3}
contains a long list of dimension combinations $(m, n)$ such there exists
a map $f \co  S^m \to S^n$ with $N (f, *) = 0$ but $\MCC (f, *) = 1$.

Clearly $\MC (f_1, f_2) \le 1$ whenever $[f_1] - [a \scirc f_2]$ lies in
$E (\pi_{m -1} (S^{n -1}))$, the image of the Freudenthal suspension. On
the other hand, it is well known that $\MC (f_1, f_2)$ is infinite if
$[f_1] - [a \scirc f_2] \not\in E (\pi_{m -1} (S^{n -1}))$ and $m, n \ne
(1, 1)$. This follows also as a very special consequence of one of our
results concerning $\MC$  (cf \fullref{3.5}(i)).
\end{example}

In \fullref{sec3} we discuss in great generality lower and upper bounds
for $\MC (f_1, f_2)$ (and deduce the classical Wecken theorem for
coincidences as a corollary). In particular, in the case $M = S^m$ we
obtain a necessary (and, if $X_m (N) = 0$ , also sufficient) finiteness
condition for $\MC (f_1, f_2)$, expressed in terms of $\omega^\# (f_1,
f_2)$ (cf \fullref{6.17} below).

\begin{example}                          
\label{1.25}
Let $N$ be an odd-dimensional spherical space form (ie the quotient
of $S^n$ by a free action of a finite group $G$). Then we have for all
$f_1, f_2 \co  S^m \to N$:
$$\MCC (f_1, f_2) = N^\# (f_1, f_2) =
\begin{cases}
0 & \text{if } f_1 \sim f_2 \text{ or } m < n; \\
\#G & \text{if } f_1 \not\sim f_2 \text{ and } m > 1; \\
|d^0(f_1)-d^0(f_2)| & \text{if } m = 1 \text{ and } N = S^1.
\end{cases}$$
(Here $d^0 (f_i) \in \mathbb Z$ denotes the usual degree).

Moreover, if $n \ge 3$ then
$$\MC (f_1, f_2) =
\begin{cases}
\infty & \text{if } [f_1] - [f_2] \not\in p_* \scirc E(\pi);\\
0 & \text{if } f_1 \sim f_2 \text{ or } m < n; \\
\#G & \text{otherwise}.
\end{cases}$$
Here
\begin{equation*}
\begin{CD}
\pi \subset \pi_{m -1} (S^{n -1}) @>{ E }>> \pi_m (S^n)
@>{ p_* }>> \pi_m (N)
\end{CD}
\end{equation*}
are the natural (eg suspension) homomorphisms and $\pi$ denotes all
of $\pi_{m -1} (S^{n -1})$ if $\# G \le 2$, and the kernel of the (total)
Hopf--Hilton homomorphism $h$ (cf \fullref{6.18} below) if $\# G \ge 3$.

Observe that no specific feature of the group action -- apart from the
order of $G$ -- enters the picture here. For a geometric explanation
of such phenomena in terms of {\em almost injective points} see the
discussion of \fullref{6.22} below.
\end{example}

For further concrete geometric settings where the minimum numbers $\MC
(f_1, f_2)$ and $\MCC (f_1, f_2)$ have been calculated explicitly, 
see Koschorke \cite[Examples~I--III]{K6}.

Having the precise finiteness criterion \fullref{6.17} for $\MC (f_1, f_2)$ at
our disposal we may ask:  what can we say in case it is satisfied?

\begin{thm}                                
\label{1.26}
Consider maps $f_1, f_2 \co  S^m \to N$ into an arbitrary $n$--manifold
such that $(m, n) \ne (2, 2)$.  Assume that $\MC (f_1, f_2)$ is finite. If
the suspension
$$E \co  \pi_{m -1} (S^{n -1}) \to \pi_m (S^n)$$
is injective (eg if $m < 2n - 2$ or $n = 2$), then
\begin{equation*}
N^\# (f_1, f_2) = \MCC (f_1, f_2) = \MC (f_1, f_2) \le \#
\pi_1 (N) .
\end{equation*}
\end{thm}

This follows from \fullref{6.24} below.

At the end of this paper we give purely homotopy theoretical descriptions
of our basic geometric coincidence invariant $\deg^\#$ (cf
\fullref{7.6} and \fullref{7.8}). We obtain a decomposition of $\deg^\#$ into
two components (cf \eqref{7.11}). One of them is the Ganea--Hopf invariant
$H_{\mathcal C}$  relative to an attaching map of top dimensional
cells in the universal covering space $\wtilde N$ of $N$. It turns
out that $H_{\mathcal C}$ is a finiteness obstruction for the minimum
number $\MC$ (and actually the only one in a dimension range depending
on the connectivity of $N$, cf \fullref{7.13} and \fullref{7.16}). Here we use homotopy
theoretical tools such as Ganea's exact EHP--sequence.

In much of our discussions we can switch freely back and forth between
base point preserving maps and homotopies and their base point free
counterparts. This is made precise in \fullref{secA}.

The approach of this paper can also be applied fruitfully to general
inverse image problems (where the submanifolds $\{ *\} \subset N$ and
$\Delta \subset N \times N$, cf \eqref{1.5}, are replaced by arbitrary closed
smooth submanifolds) or to over- and under-crossings of link maps into a
manifold of the form $N \times \mathbb R$. In the latter case we obtain
unlinking obstructions which often settle unlinking questions and which,
in addition, turn out to distinguish a great number of different link
homotopy classes (and sometimes even classify them completely). Moreover,
our approach also leads to the notion of Nielsen numbers for link maps
(cf Koschorke \cite{K4}).

\begin{notcon}
\label{1.27}
All {\em manifolds} are assumed to be Hausdorff spaces having a countable
basis; they have empty boundaries unless stated otherwise. A submanifold
$C \subset M$ with a specified trivialization of its {\em normal bundle}
$\nu (C, M)$ is called {\em framed}. In any bordism set $0$ denotes the
class represented by empty data. We will often neglect the notational
distinction between constant maps and their values. Given any topological
space $X$, $P (X)$ is the space of all paths $\theta \co  I \to X$, endowed
with the compact--open topology. $\# S$ denotes the (finite or infinite)
number of elements in a set  $S$. $E$ stands for Freudenthal suspension.
$\varphi := f^*_1 (TN) - TM \in KO (M)$ and $\wtilde\varphi := \pr^*
(\varphi) \in KO (E (f_1, f_2))$ (compare \eqref{1.10}) are the relevant virtual
coefficient bundles for our (stabilized) obstruction theory. Arbitrary
reflections on spheres are denoted by~$r$.
\end{notcon}

\section{The strong $\omega$--invariant $\omega^\# (f_1, f_2)$  and the
strong Nielsen number $N^\# (f_1, f_2)$}
\label{sec2}

Throughout this paper $f_1, f_2, f \co  M \to N$ denote  (continuous)
maps between the smooth connected (non-empty)  manifolds $M$ and $N$
without boundary, of strictly positive dimensions $m$ and $n$, resp.,
$M$ being compact.

Consider the set
\begin{equation}                                 
\label{2.1}
\Omega^\# (f_1, f_2) := \{ (C, \wtilde g, \wbar g^\#) \}
\diagup  \text{bordism} \text{ in } M \times I
\end{equation}
of bordism classes of triples of the indicated form where
\begin{enumerate}
\item[(i)]  $C$ is a closed smooth submanifold of $M$;
\item[(ii)]  $\wtilde g \co  C \to E (f_1, f_2)$ is a section of $\pr|$
(cf \eqref{1.6}), ie $\pr \scirc \wtilde g$ is the inclusion;
\item[(iii)] $\wbar g^\# \co  \nu (C, M) \cong f_1^* (TN) | C$ is a
vector bundle isomorphism which gives a (nonstabilized) description of the
normal bundle of $C$ in $M$ in terms of the tangent bundle $TN$ of $N$.
\end{enumerate}

Such triples occur very naturally when we study the coincidence
behavior of $f_1$ and $f_2$. Indeed, if the map $(f_1, f_2) \co  M \
\longrightarrow N \times N$ is smooth and transverse to the diagonal
$\Delta$ then the coincidence data \eqref{1.5}--\eqref{1.7} yield the
desired triple.

If $f_1$ and $f_2$ are arbitrary continuous maps, we apply this procedure
to a smooth map $(f'_1, f'_2)$ which approximates $(f_1, f_2)$ and is
transverse to $\Delta$. Using the techniques of \cite[Section~3]{K3}, we
see that there is a canonical bijection $\Omega^\# (f'_1, f'_2) \approx
\Omega^\# (f_1, f_2)$ induced by any sufficiently small homotopy from
$(f'_1, f'_2)$ to $(f_1, f_2)$.

In any case the resulting triple $ (C, \wtilde g, \wbar g^\#)$
determines a well-defined bordism class
\begin{equation}                                               
\label{2.2}
\omega^\# (f_1, f_2) = [C, \wtilde g, \wbar g^\#]\in \Omega^\# (f_1, f_2).
\end{equation}

The same kind of argument allows us to handle also arbitrary (possibly
\lq\lq large\rq\rq) homotopies. The result can be best expressed in the
language of functors. Consider the category $\mathfrak P$ whose objects
are continuous maps $(f_1, f_2) \co  M \to N \times N$ and whose morphisms
are equivalence classes of homotopies $F \co  (f_1, f_2) \sim (f'_1, f'_2)$;
here two homotopies $f_0, f_1$ from $(f_1, f_2)$ to $(f'_1, f'_2)$ are
called equivalent if they can be deformed continuously into one another
through such homotopies (ie at each stage of the deformation $F_t,
t \in [0, 1]$,  is a homotopy from $(f_1, f_2)$ to $(f'_1, f'_2)$).

\begin{prop}                          
\label{2.3}
The nonstabilized coincidence invariant determines a functor
$(\Omega^\#\!,\omega^\#)$ from the category $\mathfrak P$ of pairs of maps and
(deformation
classes of) homotopies to the category consisting of pointed sets and
of  bijections preserving the preferred element.
\end{prop}

If we consider the coincidence submanifold $C = C (f_1, f_2) \text{ of }
M$ just as an abstract manifold  and if we stabilize $\wbar g^\#$
to yield the {\it stable} vector bundle isomorphism
$$
\wbar g \co  T C \oplus f^*_1 (TN) | C \cong TM | C
$$
we obtain the coincidence invariant
\begin{equation}                                
\label{2.4}
\wtilde \omega (f_1, f_2) = [C, \wtilde g, \wbar g]
\in \Omega_{m-n} (E (f_1, f_2) ;\quad \wtilde\varphi := \pr^*
(f_1^* (TN) - TM))
\end{equation}
which was studied in \cite{K3}. Clearly it, too, determines a functor
as above which actually takes values in (normal bordism)  {\em groups}
with a preferred element. Stabilization yields a forgetful transformation
\begin{equation}                                       
\label{2.5}
\stab\co  (\Omega^\# (f_1, f_2), \omega^\# (f_1, f_2)) \longrightarrow
(\Omega_{m-n} (E (f_1, f_2); \wtilde\varphi), \wtilde\omega
(f_1, f_2))
\end{equation}
In the stable dimension range $m \le 2n - 2$ we are dealing
with bijections here and stabilization leads to no loss of
information. (Actually, $\wtilde\omega (f_1, f_2)$ is even the
only looseness obstruction if $m < 2n - 2\, $; cf \cite[Theorem~1.10]{K3}).
However, in general there are many situations where
the nonstabilized coincidence invariant $\omega^\# (f_1, f_2)$ turns
out to be considerably more powerful than $\wtilde\omega (f_1,
f_2)$. (This is reflected by the discussion in \fullref{1.24}). On
the other hand it is often much easier to handle the stabilized invariant
$\wtilde\omega (f_1, f_2)$: it lies in a bordism {\it group} (not
just set) and computational techniques are available (especially  for low
codimensions $(m{-}n)$, cf \cite[9.3]{K1}; compare also \cite[Section~3]{K6}).

Next, given  any bordism class $c = [C, \wtilde g, \wbar g^\#]
\in \Omega^\# (f_1, f_2)$ (cf \eqref{2.1}), let
\begin{equation}                                   
\label{2.6}
c_A = [C_A = \wtilde g^{- 1} (A), \wtilde g | C_A,
\wbar g^\# | ]
\end{equation}
denote its contribution to a given path component $A \in \pi_0 (E (f_1,
f_2))$.

\begin{definition}                             
\label{2.7}
(i)\qua We call a path component $A$ of $E (f_1, f_2)$ {\it strongly
essential} if the corresponding contribution $\omega_A^\# (f_1, f_2)$
to $\omega^\# (f_1, f_2)$ is nontrivial (ie not representable by
empty data).

(ii)\qua We define the {\it strong Nielsen number} $N^\# (f_1, f_2)$ of $f_1$
and $f_2$ to be the number of strongly essential path components $A \in
\pi_0 (E (f_1, f_2))$.
\end{definition}

This is in analogy to the (\lq\lq weak\rq\rq) Nielsen number $N (f_1,
f_2)$ which was extracted from
\begin{equation}                            
\label{2.8}
\wtilde\omega (f_1, f_2) \in \Omega_{m-n} (E (f_1, f_2);
\wtilde\varphi) \cong \bigoplus_{A \in \pi_0 (E (f_1, f_2))}
\Omega_{m-n} (A; \wtilde\varphi | A)
\end{equation}
and discussed in detail in \cite{K3}.

\begin{remark}                                
\label{2.9}
Clearly, if $\omega^\# (f_1, f_2)$ is trivial then so is $N^\# (f_1,
f_2)$. However it is conceivable that the converse  does not hold
in general; indeed, the various components  $C_A (f_1, f_2)$ of the
coincidence locus may possibly link in $M$ so that their nulbordisms
cannot be fitted together to yield {\it disjoint} embeddings into $M
\times I$.

Such complications cannot arise in the stabilized theory: $N (f_1, f_2) =
0$ if and only if $\wtilde\omega (f_1, f_2) = 0$ (much like a norm in
a vector space decides precisely whether a given vector vanishes).
\end{remark}

\begin{remark}                                 
\label{2.10}
In higher codimensions our bordism approach allows us to capture
coincidence phenomena which seem to be entirely outside the reach of the
methods of singular (co)homology theory. Already the weakened stabilized
bordism invariant
$$
\omega (f_1, f_2) := \pr_* (\wtilde\omega (f_1, f_2)) \in \Omega_{m
-n} (M; \varphi := f^*_1 (TN) - TM)
$$
(which involves neither the path space $E (f_1, f_2)$, cf \eqref{1.6}, nor
the resulting Nielsen decomposition) has lead to the solution of a
problem which corresponds to determining cohomological obstructions of
arbitrarily high order (cf the theorem in the introduction of \cite{K2}
and its corollaries).
\end{remark}

\fullref{1.11} and \fullref{1.12} of the introduction follow now from 
\eqref{1.9}, \eqref{1.13}, \eqref{2.1}, \eqref{4.7},
\eqref{5.2} and from \cite[Example~I]{K3}, or they can
be proved by refining the methods of that paper (cf also
\fullref{3.5}(iii) below).

In general the nonnegative integer $N^\# (f_1, f_2)$ contains considerably
less information than the invariant $\omega^\# (f_1, f_2)$ which, however,
has the drawback that it lies in a bordism set which varies with $f_1$
and $f_2$.

This complication can be avoided in some important settings. Given $y_0
\in N$, let $\Omega (N, y_0)^+$ denote the loop space of $N$ at $y_0$,
with an extra point $+$ added; thus $S^n \wedge (\Omega (N, y_0)^+)$
is the Thom space of the trivial $n$--plane bundle over $\Omega (N, y_0)$.

\begin{prop}[cf Hatcher--Quinn {{\cite[3.1]{HQ}}}]
\label{2.11}
Assume that $M$ is
$(m{-}n{+}1)$--con\-nec\-ted. Then any choice of points $x_0 \in M$, $y_0 \in N$
and of paths $\gamma_i$ in $N$ joining $f_i (x_0)$ to $y_0$, $i = 1,
2,$ and of a local orientation of $N$ at $y_0$ induces   a bijection
$$\Omega^\# (f_1, f_2) \approx [M, S^n \wedge (\Omega (N, y_0)^+)]$$
(involving the Pontryagin--Thom procedure) and thus allows us
to identify the coincidence invariant $\omega^\# (f_1, f_2)$ with an
element of the indicated homotopy set. These identifications commute with
the bijections induced by  homotopies of $(f_1, f_2)$, cf \fullref{2.3},
provided the paths in $N$ are chosen compatibly.

Moreover there is a canonical involution $\inv$ of the homotopy set $[M,
S^n \wedge (\Omega (N, y_0)^+)]$ such that
\begin{equation}                        
\label{2.12}
\omega^\# (f_2, f_1) = \inv (\omega^\# (f_1, f_2)) .
\end{equation}
\end{prop}

\begin{proof} Given a triple $(C, \wtilde g, \wbar g^\#)$ as in
\eqref{2.1}, a cell-by-cell argument allows us to construct a homotopy
$G\co  C \times I \to M$ from the inclusion $g = \pr \scirc \wtilde
g\co C \subset M$ to the constant map $g_1$ at $x_0$. A lift of $G$ in $E
(f_1, f_2)$, starting with $\wtilde g$, ends with a map
\begin{equation}                             
\label{2.13}
\wtilde{g}_1 \co C \longrightarrow \pr^{- 1} (\{ x_0\} ) \sim
\Omega (N, y_0) ;
\end{equation}
(compare diagram \eqref{1.6}); if $x \in C$ and $\wtilde{g} (x) = (x, \theta)$
then $\wtilde{g}_1 (x)$ can be given by the concatenated path
\begin{equation}
\label{2.14}
\begin{CD}                             
y_0  \stackrel{\gamma_1^{-1}}{\longrightarrow}
f_1 (x_0)  @>{f_1 \scirc G (x, -)^{- 1}}>>
f_1 (x)  \stackrel{\theta}{\longrightarrow}
f_2 (x) @>{f_2 \scirc G (x, -)}>>  f_2 (x_0)
\stackrel{\gamma_2}{\longrightarrow} y_0.
\end{CD}
\end{equation}
$G$ induces also a vector bundle isomorphism (cf \cite[3.1]{K3}).
\begin{equation}                              
\label{2.15}
f_1^* (TN) | C \cong g_1^* (f^*_1 (TN)) \cong C \times \mathbb R^n
\end{equation}
which we compose with $\wbar g^\#$ (cf \eqref{2.1}) to obtain a normal
framing $\wbar g^\#_1 $ of $C$ in $M$. The Pontryagin--Thom
construction transforms the resulting triple $(C, \wtilde{g}_1,
\wbar g_1^\#)$ into an element of the homotopy set $[M, S^n \wedge
(\Omega (N, y_0)^+)]$.

Our connectivity assumption guarantees that this procedure yields a
well defined bijection. It implies also that $n \ge 2$ (since $H_m (M;
\mathbb Z_2) \ne 0)$. In particular, the Thom space $S^n \wedge (\Omega
(N, y_0)^+)$ is simply connected and we may identify the elements of
$\Omega^\# (f_1, f_2)$ with base point preserving or base point free
homotopy classes, as we wish. (Note also that the choices of $x_0, \
\gamma_1$ and $\gamma_2$ do not matter in case $N$ is 1--connected).

Next we describe the involution $\inv$ of the homotopy set $[M, S^n
\wedge (\Omega (N, y_0)^+)]$ at the level of bordism classes. Given a
submanifold $C$ of $M$, a map $\wtilde{g}_1 \co C \to \Omega (N,
y_0)$ and a framing $\wbar g^\# \co \nu (C, M) \cong C \times
\mathbb R^n$, evaluate $\wtilde{g}_1$ to obtain a homotopy $h \co C
\times I \to N$ from the constant map at $y_0$ to itself. Choose a
trivialization
$\begin{CD}\wbar h^\# \co h^* (TN) @>{ \cong }>> (C \times
I) \times T_{y_0} (N)\end{CD}$ which restricts to the identity over
$C \times
\{ 0\}$ and let $\wbar h^\#_1$ denote the corresponding
automorphism of $C \times T_{y_0} (N) \cong C \times \mathbb R^n$ over
$C \times \{ 1\} = C$. Compose $\wbar g^\#_1$ with $- \wbar
h^\#_1$ and $\wtilde{g}_1$ with the involution of $\Omega (N, y_0)$
which reverses the loops. The resulting triple represents $\inv ([C,
\wtilde{g}_1, \wbar g^\#_1])$.

Now apply the whole preceding discussion to the special case where $(C,
\wtilde{g}, \wbar g^\#)$ are the coincidence data of the (generic)
pair $(f_1, f_2)$. Interchanging $f_1$ and $f_2$  clearly reverses the
path $\wtilde{g}_1 (x)$, $x \in C$ (cf \eqref{2.14} where $\theta = \const$
in this case). Also if we describe $\nu (C, M)$ by $f_2^* (TN)$ (instead
of $f^*_1 (TN)$,  cf \eqref{2.1}(iii)) we must base the analogue of
\eqref{2.15} on
a trivialization $(f_2 \scirc G)^* (TN) \cong C \times I \times T_{y_0}
(N)$. Moreover note that the two projections from the diagonal $\Delta
\subset N \times N$ to $N$ yield isomorphisms $\nu (\Delta, N \times N)
\cong TN$ which play a role in \eqref{1.7} but differ by a factor $-1$. Indeed,
for all $y \in N$ and $v  \in T_y (N)$ the two vectors $(v, 0), (0, -
v) \in T_{(y,y)} (N \times N)$ yield the same element in $\nu (\Delta,
N \times N)$ since their difference is tangential to $\Delta$. This
explains the negative sign in the definition of the involution $\inv$.
\end{proof}

Our chosen local orientation of $N$ at the point $y_0 \in N$ (cf
\fullref{2.11})
determines a collapsing map
\begin{equation}                              
\label{2.16}
\coll \co N \longrightarrow N / (N - \mathring{B}{^n})
\cong B^n / \partial B^n \cong S^n
\end{equation}
up to homotopy (where $B^n$ is a small ball around $y_0$). Consider also
the maps
$$
\begin{CD}
S^n @>{ \In_\ell }>> S^n \wedge (\Omega (N, y_0)^+) @>{ \ret' }>> S^n
\end{CD}
$$
defined by the inclusion at a fixed loop $\ell \in \Omega (N, y_0)$
and by the \lq\lq horizontal projection\rq\rq. We obtain the diagram of
(induced) maps
\begin{equation}                            
\label{2.17}
  \xymatrix{
     [(M, x_0), (N, y_0)] \ar[r]^-{\deg^\#}
     \ar@<-2pt> `d[r] `[rr]_-{\coll_*} [rr] &
     [M,S^n\wedge(\Omega(N,y_0))^+]
     \ar@<0.5ex>[r]^-{\ret'_*} &
     [M,S^n]
     \ar@<0.5ex>[l]^-{\In_{\ell*}}
  }
\end{equation}

\begin{prop}                                
\label{2.18}
Assume that $M$ is $(m{-}n{+}1)$--connected.

Then we have in diagram \eqref{2.17}
$$
\begin{aligned}
\coll_* =& u_{\varepsilon *} \scirc \ret'_* \scirc \inv \scirc \deg^\#
\qquad \text{and} \\
\ret'_* \scirc \In_{\ell *} =& \text{identity}
\end{aligned}
$$
where $u_\varepsilon$ denotes a selfmap of $S^n$ of degree $\varepsilon =
(- 1)^n$ (cf also \fullref{2.11}).

Hence for every map $f \co (M, x_0) \to (N, y_0)$ the degree $\deg^\# (f)
= \omega^\# (f, y_0)$ determines the homotopy class of $\coll \scirc f\co 
M \to S^n$. Moreover the selfcoincidence invariant $\omega^\# (f, f)$
lies in the image of $\In_{\ell *}$ for some $\ell \in \Omega (N, y_0)$
and therefore is determined by the seemingly weaker invariant $\ret'_*
(\omega^\# (f, f))$.

In particular, if $N = S^n$ then  $\deg^\# (f)$ determines $f$ up to
homotopy and we have
$$N^\# (f, y_0) = \MCC (f, y_0) =
\begin{cases}
0 & \text{if } f \sim y_0; \\
1 & \text{if } f \not\sim y_0.
\end{cases}$$
\end{prop}

\begin{proof} Generically $f$ is smooth with regular value
$y_0$. Represent $\omega^\# (y_0, f) = \inv (\deg^\# (f))$ (cf
\eqref{2.12})
by the triple $(C = C (y_0, f) = f^{- 1} (\{ y_0\} ), \wtilde{g}_1,
\wbar g^\#_1)$ (cf \eqref{2.13}, \eqref{2.15}) and forget the map
$\wtilde{g}_1 \co C \to \Omega (N, y_0)$ (or, equivalently, apply the \lq\lq
horizontal\rq\rq\ projection $S^n \wedge (\Omega (N, y_0)^+) \to
S^n)$. Then the data $(C, - \wbar g^\#_1)$ give the standard
characterization (in the sense of Pontryagin--Thom) of $\coll \scirc
f$. (For the negative sign compare the proof of \fullref{2.11}.)

An inspection of \eqref{2.14} shows that $\wtilde{g}_1$ is homotopy trivial
when $f_1 = f_2$ and $\theta$ is a constant path. Thus $\omega^\# (f, f)$
lies in the image of $\In_{\ell *}$ where $\ell = \gamma^{- 1}_1 \gamma_2$
(compare \fullref{2.11}).

If $N = S^n$, then $\# \pi_0 (E (f, *)) = 1$ must agree with $N^\# (f,
*)$ whenever $f \not\sim *$ and hence $\deg^\# (f)$ is nontrivial. The
last claim follows from \fullref{1.11}(iii) and can be established also for $n =
2$ by a connected sum argument.
\end{proof}

\section{Isolated coincidence points}                    
\label{sec3}

In this section we establish a finiteness criterion as well as upper
bounds for $\MC (f_1, f_2)$. In the special case $m = n > 2$ we reobtain
the classical Wecken theorem.

In view of \fullref{1.12} we may (and do) assume that $m \ge 2$.

First we construct a map
\begin{equation}                                 
\label{3.1}
e \co \bigoplus_{A \in \pi_0 (E (f_1, f_2))} \pi_{m -1} (S^{n -1})
\longrightarrow \Omega^\# (f_1, f_2)
\end{equation}

\begin{stepI}  Given a path component $A$ of $E (f_1, f_2)$, choose an
element $(x, \theta) \in A$ as well as the following data:

(i)\qua a diffeomorphism $b \co  (B^m, 0) \to (B_x, x)$ from the standard
compact unit ball $B^m$ in $\mathbb R^m$ onto a small \lq\lq ball\rq\rq\
$B_x$ around $x$ in $M$; and

(ii)\qua a trivialization $\wbar b$  of $f^*_1 (TN) | B_x$, the pullback
of the tangent bundle of $N$ by $f_1 | B_x$.

These data, together with the inclusion $B_x \subset M$,
allow us to interpret any framed closed smooth submanifold $C$ of
$\mathring{B}^m$ as a triple $(C, g, \wbar g^\#)$ (compare
\eqref{1.6},
\eqref{1.7}, and \eqref{2.1}). Moreover the embedding $g \co C \subset M$ lifts canonically
(up to homotopy) to a map $\wtilde{g}$ into $E (f_1, f_2)$ as follows:
given $x' \in C$, let $r_{x'}$ be the straight path in $B_x \approx B^m$
joining $x'$ to $x$ and define $\wtilde{g} (x') = (x', (f_1 \scirc
r_{x'}) \theta (f_2 \scirc r_{x'})^{- 1})$ by concatenating the resulting
image paths in $N$ with the fixed chosen path $\theta$. In view of the
Pontryagin--Thom procedure we obtain the well-defined composite map
\begin{equation}                                 
\label{3.2}
\begin{CD}
e_A \co  \pi_{m -1} (S^{n -1}) @>{E}>> \pi_m (S^n)
\longrightarrow \Omega^\# (f_1, f_2)
\end{CD}
\end{equation}
where $E$ is the Freudenthal suspension.

Since tubular neighborhoods are essentially unique up to isotopy, the
map $e_A$ does not really depend on the data $b$ and $\wbar b$,
but only on the resulting orientations of the tangent spaces $T_x M$
and $ T_{f_1 (x)} N$. Clearly, changing {\it one} of these
orientations amounts to replacing $e_A (\alpha)$ by $e_A (- \alpha), \
\alpha \in \pi_{m -1} (S^{n -1})$.Thus, in fact, $e_A$ depends only
on the co-orientation of $f_1$ at $x$, ie on the orientation of
the virtual coefficient bundle $\varphi = f^*_1 (TN) - TM$ (cf
\fullref{1.27})
at $x$ (or, equivalently, of $\wtilde\varphi = \pr^* (\varphi)$ at
$(x, \theta)$, cf \eqref{2.4}). Furthermore, $e_A$ does not really depend
on the choice of the element $(x, \theta)$ in $A$ either. Indeed, any
path in $A$ (if it projects to a smooth path in $M$) gives rise to an
isotopy of $B_x$ and to the corresponding deformation of the other
data needed in our construction.
\end{stepI}

\begin{stepII}Given any element $\{ \alpha _{A_j}\} \in \oplus \pi_{m -1}
(S^{n -1})$ in the domain of $e$ (cf \eqref{3.1}), the summands $\alpha_{A_j}$
are nontrivial for only finitely many path components $A_1, \dots, ,
A_k$ of $E (f_1, f_2)$. Apply the construction in step I to each of them
simultaneously, using disjoint balls $B_{x_1}, \dots, B_{x_k}$, and define
\begin{equation}                                
\label{3.3}
e (\{ \alpha_{A_j} \} ) = \coprod^k_{j = 1} e_{A_j} (\alpha_{A_j})
\in \Omega^\# (f_1, f_2).
\end{equation}
If we can (and do!) orient the virtual coefficient bundle
$\wtilde\varphi$ (cf \eqref{2.4}) then the map $e$ is well defined and
independent of all other choices (recall our assumption  $m \ge 2$
which allows for the necessary disjoint isotopies). In general $e$ is
only well defined up to replacing some of the summands $\alpha_{A_j}$
by  $- \alpha_{A_j}$. But in any case the image of $e$ is a well defined
subset of $\Omega^\# (f_1, f_2)$.

Furthermore the construction of $e$ and of its image is compatible with
homotopies of $(f_1, f_2)$ (cf \fullref{2.3}).
\end{stepII}

Now consider an {\it isolated} coincidence point $x \in M$ of $f_1$
and $f_2$. Identify a neighborhood $U_y$ of $y := f_1 (x) = f_2 (x)$
in $N$ with the Euclidean space via a diffeomorphism $(U_y, y) \approx
(\mathbb R^n, 0)$; also identify $B^m$ with a small ball $B_x \subset
f^{-1}_1 (U_y) \cap f^{-1}_2 (U_y)$ around $x$ (as in step I, (i), above)
which contains no coincidence point other than $x$. This allows us to
define the \lq\lq index map\rq\rq\
\begin{equation}                                 
\label{3.4}
\begin{CD}
q = \frac{f_1 - f_2}{\Vert f_1 - f_2\Vert} \co  S^{m -1}
@>{\cong}>{b|}> \partial B_x \longrightarrow S^{n -1}
\end{CD}
\end{equation}
(compare Wyler \cite[Theorems~3 and~4]{Wy}) and its \lq\lq concentric
extension\rq\rq\ $Q \co  B_x \to B^n$, $Q (t b(z)) := t q(z)$ for $z
\in S^{m -1}$. We can easily deform $f_1$ (while leaving it fixed
outside of a small neighborhood of $B_x$) into a map $f'_1$ such
that the corresponding index map $q'$ is smooth and $f'_1 = f_2 + Q'$
on $B_{x}$. Indeed, leaving the map $f_2$ (which we may assume to be
smooth) unchanged, first deform the difference $f_1 - f_2$ in a collar
neighborhood of $S^{m -1}$ until $(f_1 - f_2) | S^{m -1}$ coincides with
a smooth approximation $q'$ of $q$; then use the linear structure on
$U_y \approx \mathbb R^n$ to obtain the desired homotopy from $(f_1 -
f_2) | B_x$ to  the concentric extension $Q'$ of $q'$.
After a further small homotopy near $x$ we have $f'_1 = f_2 + Q' -
\varepsilon \cdot *$ where $* \in S^{n -1}$ is a regular value of
$q'$. Thus the (generic!) coincidence locus of $f'_1$ and $f_2$ near
$x$ consists of the framed submanifold $\varepsilon q^{' -1} (\{*\})$
of the $\varepsilon$--sphere around $x$ and hence corresponds to the
Freudenthal suspension of the homotopy class of $q' \sim q$.

We conclude that the contribution of an isolated coincidence point $x \in
C (f_1, f_2)$ towards $\omega^\# (f_1, f_2)$ is given by the image $e_A
([q])$ of its index map $q$ (cf \eqref{3.2} and \eqref{3.4}; here $A$ contains $(x$,
constant path at $f_1 (x) = f_2 (x))$).

\begin{thm}                              
\label{3.5}
 Assume that the minimum number $\MC (f_1, f_2)$ is finite. Then the
 following holds.

{\rm (i)}\qua If $m \ge 2$, then $\omega^\# (f_1, f_2)$ lies in the image
of the map $e$ (cf \eqref{3.1}). In particular, if $m > n$ the homology
class defined by a generic coincidence manifold $C (f_1, f_2)$ in $M$
must vanish.

{\rm (ii)}\qua If $n \ne 2$, then $\MC (f_1, f_2) \le \# \pi_0 (E (f_1,
f_2))$; more precisely: if $n > 2$ and $C (f_1, f_2)$ is finite, then
there exists a deformation of $f_1$ and $f_2$ which preserves empty
Nielsen classes of coincidence points and replaces each nonempty Nielsen
class by (at most) a single coincidence point.

{\rm (iii)}\qua If $\pi_{m -1} (S^{n -1}) = 0$ (eg if $m > n = 1$
or $2$) then $\MC (f_1, f_2) = 0$, ie $(f_1, f_2)$ is homotopic to
a coincidence free pair of maps.
\end{thm}

\begin{proof} We may assume that $C (f_1, f_2)$ is finite. Let $x, x_1,
\dots, x_r$ be the coincidence points whose $\wtilde{g}$--values (cf
\eqref{1.6}) lie in a given path component $A$ of $E (f_1, f_2)$.

Recall that the construction of $\omega^\# (f_1, f_2)$ involves a generic
approximation of $(f_1, f_2)$. As a result each isolated coincidence
point $x_j$ is replaced by an $n$--codimensional coincidence submanifold
$C_j$ (eg of the form $C_j = \varepsilon q^{' -1}_j (\{*\})$ as in
the discussion following \eqref{3.4}) which lies in a small neighborhood $U_j
\approx \mathbb R^m$ of $x_j$ in $M$, $j = 1, \dots, r$.
In order to prove claim (i) we have to \lq\lq slide\rq\rq\ these
(disjoint) neighborhoods into a small ball $B_x \subset M$ around
$x$. For this purpose pick a smooth \lq\lq Nielsen\rq\rq\ path $\sigma_j$
from $x_j$ to some point $x'_j \in B_x$ such that $\sigma_j$ avoids all
coincidence points except $x_j$. (Here the Nielsen property of $\sigma_j$
means that a lifting to $E (f_1, f_2)$, cf \eqref{1.6}, joins the constant
path at $f_1 (x_j) = f_2 (x_j)$ to a short path near $f_1 (x) = f_2
(x)$). The corresponding arc
$$
a := \{ (\sigma_j (t), t) | t \in I \}
$$
has a tubular neighborhood of the form $a \times \mathbb R^m$ in $M
\times I$. Then the resulting submanifold $a \times C_j$ gives rise to
a bordism which relates the local $\omega^\#$--data of $x_j$ (ie
the triple $(C_j, \wtilde{g}_j, \wbar g^\#_j)$ of (partial)
coincidence data corresponding to $C_j$, cf \eqref{1.5}--\eqref{1.7})
to a contribution of the form $e_A ([q_j])$ towards $\omega^\# (f_1, f_2)$
(compare \eqref{3.2} and \eqref{3.4}). Since $E$ preserves addition we may iterate our
argument to show that all of $\omega^\# (f_1, f_2)$ lies in the image
if $e$.

In order to prove claim (ii) we may have to deform $f_1$ and $f_2$
(and not just the local $\omega^\#$--data). In view of \fullref{1.12}
we need to consider only the case when $m, n \ge 3$. Then, given two
isolated coincidence points $x, x_j$ in the same Nielsen class, we can
join them by a smoothly embedded Nielsen path $\sigma$ which misses
all other coincidence points.  Similarly $f_1 \scirc \sigma$ will be an
embedding after a small deformation of $f_1$ in a tubular neighborhood
of $\sigma (I)$. Moreover, by making an approximation transverse to
$f_1 \scirc \sigma (I)$ and possibly by shifting $I$--levels slightly
near finitely many intersection points we obtain a  homotopy $F$ from
$f_1 \scirc \sigma$ to $f_2 \scirc \sigma$ which leaves endpoints fixed
and which coincides with the constant homotopy of $f_1 \scirc \sigma$
only at the initial deformation parameter.

Next identify a small tubular neighbourhood $V \subset M$ of the arc
$\sigma (I)$ with $I \times \mathbb R^{m -1}$. Without changing anything
outside of $V$ we want to use the homotopy $F$ in order to deform $f_2$
into a map $f'_2$ such that $f_1$ and  $f'_2$ coincide in $\sigma (I) =
I \times \{ 0\}$ but nowhere else in $V$. Pick a smooth map $\rho \co I \to
[0, \infty)$ such that $\rho^{- 1} (\{ 0\} ) = \{ 0, 1\}$. First deform
the identity map of $M$ along the normal rays in $V$ towards $\sigma
(I)$ until each point $(t, v)$ of
$$
V_\rho := \{ (t, v) \in V = I \times \mathbb R^{m -1} | \Vert v
\Vert \le \rho (t) \}
$$
gets mapped to $(t, 0)$. Compose this homotopy with $f_i$, $i = 1,
2$, so that in the end  $f_i (t, v) = f_i \scirc \sigma (t)$ whenever
$\Vert v \Vert \le \rho (t)$. Now define $f'_2$ to equal $f_2$ outside
of $V_\rho$, and $f'_2 (t, v) := F \bigl(t, \frac{\Vert v\Vert}{\rho
(t)}\bigr)$
when $(t, v) \in V_\rho, t \ne 0, 1$ (ie along any normal ray in
$V_\rho$, starting at $\sigma (t)$, the constant path with value $f_2
(\sigma (t))$ is replaced by the path $F (t, -)$ from $f_1 (\sigma (t))$
to $f_2 (\sigma (t)))$.

Clearly $f'_2$ is homotopic to $f_2$. We have replaced our two original
coincidence points of $(f_1, f_2)$ by the full arc $\sigma (I)$ of
coincidence points of $(f_1, f'_2)$. Extending $\sigma$ to a slightly
larger interval $I_+ = (- \varepsilon, 1 + \varepsilon)$ and using
suitable tubular neighbourhoods of $\sigma (I_+)$ and $f_1 \scirc
\sigma (I_+)$ we can find sets $B_x \subset M$ and $U_y \subset N$
as in the discussion of \eqref{3.4} such that $\sigma (I) \subset B_x$. After
further (local) homotopies $f_1 - f'_2$ is the concentric extension $Q$
of some index map $q$ (cf \eqref{3.4}). Thus in the end $x$ and $x_j$ are
replaced by just one coincidence point at the center of $B_x$. Iterating
this procedure we can reduce each Nielsen class to a single point (see
\cite[Section~3]{K3} for the compatibilities of Nielsen decompositions
with homotopies of $f_1$ and $f_2$).

Finally assume that $\pi_{m -1} (S^{n -1})$ is trivial. Then at each
isolated coincidence point $x$ the index map (cf \eqref{3.4}) allows an
extension $Q'$ without zero. The resulting map $f'_1 = f_2 + Q'$ is
homotopic to $f_1$ and has no coincidence with $f_2$.
\end{proof}

\begin{example}[$ m = n \ge 2$]
\label{3.6}
Here we may identify $\pi_{m -1} (S^{n
-1})$ with $\mathbb Z$ via the mapping degree. Then the map
$$e \co  \bigoplus\limits_{A\in\pi_0(E(f_1,f_2))}
  \xymatrix{
     \mathbb{Z} \ar[r] &
    \Omega^\#(f_1,f_2) \ar@{}[r]|-{=} &
    \Omega_0(E(f_1,f_2);\wtilde{\varphi}) \ar@{}[r]|-{=} &
    \bigoplus\mathbb{Z}\oplus\bigoplus\mathbb{Z}_2  \\
    & \omega^{\#}(f_1,f_2) \ar@{}[r]|-{=}
    \ar@{}[u]|-{\rotatebox{90}{$\in$}} &
    \tilde{\omega}(f_1,f_2) \ar@{}[u]|-{\rotatebox{90}{$\in$}}
  }
$$
is the direct sum of identity maps and  $\mod 2$ reductions according
to the orientability of the coefficient bundle $\wtilde\varphi|A$
(compare the discussion following \eqref{3.2}, and \cite[9.3]{K1}).

Now assume that $n > 2$. After suitable deformations the Nielsen class
$C_A = \wtilde{g}^{- 1} (A)$ corresponding to a path component $A$
of $E (f_1, f_2)$ consists of at most one point $x_A$ (cf \fullref{3.5}(ii)).

If $A$ is nonessential we may remove the coincidence at $x_A$
altogether. This is clear when $\wtilde\varphi | A$ is oriented since
then the index map of $x_A$ has degree $0$ and hence is nulhomotopic. If
$\wtilde\varphi | A$ is not orientable this degree is even; therefore
we may replace $x_A$ by nearby generic coincidences which occur in pairs
of points having the same sign $+ 1$ or $- 1$; join each such pair by an
embedded Nielsen path $\sigma$ which reverses the given local orientation
of $\wtilde\varphi | A$; in a tubular neighbourhood of $\sigma (I)$
the two endpoints have opposite signs and give rise to a nulhomotopic
index map (compare the proof of \fullref{3.5}(ii)).

We conclude that $\MC (f_1, f_2) \le N (f_1, f_2)$. In view of
\fullref{1.11}(iii) we obtain the following classical \lq\lq Wecken
theorem\rq\rq\ as a  special consequence of \fullref{3.5}.
\end{example}

\begin{cor}                                               
\label{3.7}
If $m = n \ne 2$, then for all maps $f_1, f_2 \co  M^m \to N^n$ we have
$$
N (f_1, f_2) = N^\# (f_1, f_2) = \MCC (f_1, f_2) = \MC (f_1,
f_2) .
$$
\end{cor}

For $m = n = 1$ this follows from \fullref{1.12}.

\begin{remark}                                    
\label{3.8}
When $m \ge 1$ and $n = 1$ or $2$ the minimum number $\MC (f_1, f_2)$ --
if it is known to be finite -- is rather easily determined (cf
\fullref{1.12} and \fullref{3.5}(iii)) {\it except} in the dimension setting $m =
n = 2$ originally studied by J Nielsen. Here it took 57 years until
the central question \lq\lq are $N (f_1, \id)$ and $\MC (f_1, \id)$
always equal?\rq\rq\ was proved by B Jiang \cite{Ji1,Ji2} to have a
negative answer. Thus when $m = n = 2$ it is {\it not}
always possible to replace each Nielsen class by at most one coincidence
point (compare \fullref{3.5}(ii), and \fullref{3.6}).
\end{remark}

\section{The root case}                      
\label{sec4}

Let a map $f \co M \to N$ be given. In this section we discuss the
Nielsen and minimum numbers of a pair of the form ($f_1 = f, f_2 =$
constant map).

Since $N$ is path-connected these numbers are independent of the constant
value of $f_2$. Thus, given a basepoint $x_0 \in M$, put $y_0 := f (x_0)$;
we may assume that $f_2 = y_0$ (our notation will not distinguish between
a constant map and its value).

We define
\begin{equation}                         
\label{4.1}
\begin{aligned}
\deg^\# (f)& := \omega^\# (f, y_0) \in \Omega^\# (f, y_0), \\
\widetilde\deg (f)& := \wtilde\omega (f, y_0) \in \Omega_{m-n}
(E (f, y_0); \wtilde\varphi), \text{ and} \\
\deg(f)& := \omega (f, y_0) \in \Omega_{m-n} (M; \varphi :=
f^*_1 (TN) - TM)
\end{aligned}
\end{equation}
(compare \eqref{1.13}, \fullref{2.10} and \eqref{2.4})

According to \cite[2.1]{K3} we have
\begin{equation}                       
\label{4.2}
\# \pi_0 (E (f, y_0)) = [\pi_1 (N, y_0) : f_* (\pi_1 (M, x_0))].
\end{equation}
ie the cardinality of $\pi_0 (E (f, y_0))$ or, equivalently, of the
Reidemeister set $R (f, y_0)$, equals the index of the subgroup $f_*
(\pi_1 (M, x_0))$ in $\pi_1 (N, y_0)$. Our analysis will be based on the
simple  but useful observation that the path components of $E (f, y_0)$
can in fact be parametrized by the endpoints of certain liftings in an
appropriate covering space of $N$ (cf \fullref{4.5} below).

Consider the commuting diagram
\begin{equation}                          
\label{4.3}
  \xymatrix{
    & (\wtilde{N}_f,\wtilde{y}_0) \ar[d]^-p \\
    (M,x_0) \ar[ur]^-{\wtilde{f}} \ar[r]_-{f} &
    (N,y_0)
  }
  \end{equation}
where $p$ denotes the (basically unique) covering such that
$$p_* (\pi_1 (\wtilde{N}_f, \wtilde{y}_0))  = f_* (\pi_1 (M, x_0))$$
(cf Greenberg--Harper \cite[6.9]{GH}) and $\wtilde{f}$ is the lifting of
$f$ determined by the choice of $\wtilde{y}_0 \in p^{- 1} (\{ y_0\} )$.

\begin{lem}                           
\label{4.4}
For every point $\wtilde{y} \in \wtilde{N}_f$ the space $E (\wtilde
f, \wtilde{y})$ is path-connected.
\end{lem}

\begin{proof} Since the projection from $E (\wtilde{f}, \wtilde{y})$
to $M$ is a Hurewicz fibration we have only to join two elements of
the form $(x_0, \wtilde\theta_i)$ in $E (\wtilde{f}, \wtilde
y)$, $i = 1, 2$. Now the path $\wtilde\theta_1$ is homotopic $\rel
(0,1)$ to a concatenated path $\wtilde\tau \wtilde\theta_2$ where
$\wtilde\tau$ starts and ends at $\wtilde{y}_0$. Thus there is a loop
$\sigma$ in $M$ starting and ending at $x_0$ such that $\wtilde\tau
\sim \wtilde{f} \scirc \sigma$. Clearly $\sigma$ lifts to a path in $E
(\wtilde{f}, \wtilde{y})$ which joins $(x_0, \wtilde\theta_1)$
to $(x_0, \wtilde\theta_2)$.
\end{proof}

The preceding lemma allows us to study our original map $f$ with the
help of $\wtilde{f}$. Given $(x, \theta) \in E (f, y_0)$, define
$\psi (x, \theta) = \wtilde\theta (1)$ where $\wtilde\theta \in
P (\wtilde{N}_f)$ is the lifting of $\theta \in P (N)$ (cf \eqref{1.6})
satisfying $\wtilde\theta (0) = \wtilde{f} (x)$.

\begin{cor}                 
\label{4.5}
$\psi$ induces a bijection
$$
\pi_0 (E (f, y_0)) \longrightarrow p^{- 1} (\{ y_0\} )
$$
which is compatible with homotopies of $f$ (compare \cite[3.2]{K3}).
\end{cor}

In fact the lifting procedure $\theta \longleftrightarrow
\wtilde\theta$ determined by $\wtilde{f}$ yields a homeomorphism
\begin{equation}                                   
\label{4.6}
E (f, y_0) \cong \coprod_{\wtilde{y} \in p^{-1} (\{ y_0\} )}
E (\wtilde{f}, \wtilde{y}).
\end{equation}
Note that all these path components are homotopy equivalent via
concatenation: any path $\wtilde\tau$ in $\wtilde{N}_f$ from
$\wtilde{y}_1$ to $\wtilde{y}_2$ induces a homotopy equivalence
\begin{equation}                               
\label{4.7}
\wtilde\tau_\# \co  E (\wtilde{f}, \wtilde{y}_1) \we
E (\wtilde{f}, \wtilde{y}_2), \qquad  \wtilde\tau_\#
(x, \wtilde\theta ) = (x, \wtilde\theta \wtilde\tau) .
\end{equation}

\begin{thm}                  
\label{4.8}
Whenever $n \ge 1$ the following statements are equivalent:
\begin{enumerate}
\item[{\rm (i)}] $\deg^\# (f) = 0$;
\item[{\rm (ii)}] $N^\# (f, y_0) = 0$;
\item[{\rm (iii)}] $N^\# (f, y_0) \ne \# \pi_0 (E (f, y_0))$;
\item[{\rm (iv)}] $\deg^\# (\wtilde{f}) = 0$;
\item[{\rm (v)}] $N^\# (\wtilde{f}, \wtilde{y}_0) = 0$ .
\end{enumerate}

In particular, if $\pi_0 (E (f, y_0))$ is infinite or if the manifold
$N$ (or $\wtilde{N}_f$) is noncompact, then $\deg^\# (f), \deg^\#
(\wtilde{f})$ and the Nielsen numbers $N^\# (f, y_0)$ and $N^\#
(\wtilde{f}, \wtilde{y}_0)$ vanish.
\end{thm}

\begin{remark}                 
\label{4.9}
Precisely the analogous statement holds for the (\lq\lq stabilized\rq\rq
) degrees and Nielsen numbers $\widetilde\deg (f), \widetilde\deg
(\wtilde{f}), N (f, y_0)$ and $N (\wtilde{f}, \wtilde{y}_0)$
(cf \cite[1.11]{K3}).

In the nonstabilized setting the implication (ii) $\Longrightarrow$
(i) is nontrivial (see \fullref{2.9}).
\end{remark}

\begin{proof}[Proof of \fullref{4.8}] Clearly (i) $\Longrightarrow$
(ii) $\Longrightarrow$ (iii) and (iv) $\Longleftrightarrow$ (v). Also
the full claim of the theorem holds for $N \cong S^1$ (cf \fullref{1.12})
and (trivially) for $N \cong \mathbb R$.

Thus assume that $n \ge 2$.

In any element $[C, \wtilde{g}, \wbar g^\#]$ of $\Omega^\# (f,
y_0)$ or $\Omega^\# (\wtilde{f}, \wtilde{y}_0)$ the homotopy defined
by $\wtilde{g}$  induces a vector bundle isomorphism
\begin{equation}                     
\label{4.10}
f^* (TN) | C \cong \const_{y_0}^* (TN) | C = C \times V
\end{equation}
where $V := T_{y_0} (N)$. Thus we may henceforth interpret $\wbar
g^\#$ as a trivialization of the normal bundle $\nu (C, M)$ (compare
\eqref{1.7}). This gives rise to an identification of a tubular neighborhood $U$
of $C$ in $M$ with $C \times V$.

Now, given path components $A_1, \dots, A_k$ and $A$ of $E (f, y_0)$,
we construct maps
\begin{equation}                             
\label{4.11}
\begin{CD}
\Omega^\# (\wtilde{f}, \wtilde{y}_0) @>{ \pinch_{\mathfrak A} }>>
\Omega^\# (f, y_0) @>{ \forg_A }>> \Omega^\# (\wtilde{f},
\wtilde{y}_0)
\end{CD}
\end{equation}
as follows.

Let $\{ \wtilde{y}_1, \dots, \wtilde{y}_k\} \subset p^{- 1}
(\{ y_0\} )$ correspond to $\mathfrak A := \{ A_1, \dots, A_k\}$ via
$\psi$ (cf \fullref{4.5}). Choose pairwise distinct points $z_i$
in $V$ as well as paths $\wtilde\tau_i$ in $\wtilde{N}_f$ from
$\wtilde{y}_0$ to $\wtilde{y}_i$, $i = 1, \dots, k$. For $c = [C,
\wtilde{g}, \wbar g^\#] \in \Omega^\# (\wtilde{f}, \wtilde
y_0)$ we represent $\pinch_{\mathfrak A} (c)$ by the union of the \lq\lq
parallel\rq\rq\ submanifolds
\begin{equation}                           
\label{4.12}
C_i := C \times \{ z_i\} \subset U
\end{equation}
of $M$, together with the data $\wtilde{G}$ and $\wbar G$ described
as follows (compare \eqref{1.6} and \eqref{1.7}). In order to obtain $\wtilde{G} (x,
z_i)$, $x \in C$, apply the covering map $p$ to the concatenated path
$$\xymatrix{
      \wtilde{f}(x,z_i) \ar@{~>}[r]^-{\wtilde{f}\circ\ell} &
      \wtilde{f}(x) \ar@{~>}[r]^-{\wtilde{g}(x)} &
      \wtilde{y}_0 \ar@{~>}[r]^-{\wtilde{\tau}_i} &
      \wtilde{y}_i}$$
in $\wtilde{N}_f$; here $\ell$ denotes the linear path in $\{ x\}
\times V$ from $(x, z_i)$ to $x = (x, 0)$. $\wbar G$ is the obvious
trivialization of the normal bundle of $C_i$ in $C \times V = U \subset
M$, composed with a fixed reflection in $V$ if $TN$ is not orientable
along the closed loop $p \scirc \wtilde\tau_i$ in $N$.

Similarly, choose a path $\wtilde\tau$ in $\wtilde{N}_f$ from
$\wtilde{y}_0$ to $\wtilde{y} := \psi (A) \in p^{-1} (\{
y_o\} )$ (cf \fullref{4.5}). The map $\forg_A$ in \eqref{4.11} is then built up in three
steps. First it forgets all but the  $A$--component
\begin{equation}                                 
\label{4.13}
\omega^\#_A = [ \wtilde{g}^{- 1} (A), \wtilde{g}|, \wbar
g^\# | ]
\end{equation}
of an element $\omega^\# = [C, \wtilde{g}, \wbar g^\#] \in \Omega^\#
(f, y_0)$. Lifting $\wtilde{g}|$ to $\wtilde{N}_f$ then makes
$\omega_A^\#$ into an element of $\Omega^\# (\wtilde{f}, \wtilde{y})$
(cf \eqref{4.6}). Finally concatenate with $\wtilde\tau^{- 1}$  to obtain
$\forg_A (\omega^\#) \in \Omega^\# (\wtilde{f}, \wtilde{y}_0)$.

If $A = A_i$ for some $1 \le i \le k$ and if the chosen paths
$\wtilde\tau$ and $\wtilde\tau_i$ are homotopic $\rel (0, 1)$, then
\begin{equation}
\label{4.14a}
\forg_A \scirc \pinch_{\mathfrak A} = \text{identity}.
\end{equation}
However, if $A \not\in \mathfrak A = \{ A_1, \dots, A_k\}$ then
\begin{equation}
\label{4.14b}
\forg_A \scirc \pinch_{\mathfrak A} \equiv 0.
\end{equation}
\stepcounter{equation}
Let us apply this whole discussion to $\deg^\# (f) = \omega^\# (f,
y_0)$ and $\deg^\# (\wtilde{f}) = \omega^\# (\wtilde{f}, \wtilde
y_0)$. Identify a small neighborhood $U_{\wtilde{y}_0}$ of $\wtilde
y_0$ in $\wtilde{N}_f$ with $V$. After a suitable homotopy we may
assume that $\wtilde{f}$ is smooth, with regular value $\wtilde
y_0$, and agrees on a tubular neighborhood $U \cong C \times V$ of $C :=
\wtilde{f}^{-1} (\{ \wtilde{y}_0\} )$ with the projection to $V =
U_{\wtilde{y}_0}$.

In the construction of $\pinch_{\mathfrak A}$ (cf \eqref{4.11}) choose the set
$\mathfrak A = \{ A_1, \dots, A_k\} \subset$\break $\pi_0 (E (f, y_0))$ to be so big
that $\wtilde{f}$ avoids $F^- := p^{- 1} (\{ y_0\} ) - \{ \wtilde
y_1, \dots, \wtilde{y}_k\}$. Since $n \ge 2$ we may also choose the
paths $\wtilde\tau_i$ to go straight from $\wtilde{y}_0$ to $z_i$
in $U_{\wtilde{y}_o} = V$, and then to $\wtilde{y}_i$ via disjoint
embedded arcs in $\wtilde{N}_f - F^-$, $i = 1, \dots, k$.

Now compose $\wtilde{f}$ with an isotopy of $\wtilde{N}_f$ which
leaves $F^-$ fixed and which moves $z_i$ along $\wtilde \tau_i$
to $\wtilde{y}_i$. At the final stage of this deformation the
$\deg^\#$--invariants of $\wtilde{f}$ and $f = p \scirc \wtilde{f}$
satisfy the relation
\begin{equation}                              
\label{4.15}
\deg^\# (f) = \pinch_{\mathfrak A} (\deg^\# (\wtilde{f}))
\end{equation}
(cf \eqref{4.12}). Thus (iv) $\Rightarrow$ (i) in \fullref{4.8}. On the other hand: if
$\deg^\#_A (f)$ is trivial for some path component $A$ of $E (f, y_0)$
(cf \eqref{2.6})  then so is
\begin{equation}                                       
\label{4.16}
\deg^\# (\wtilde{f}) = \forg_A (\deg^\# (f))
\end{equation}
(cf \eqref{4.11}). Hence (iii) $\Rightarrow$ (iv) and \fullref{4.8} follows.
\end{proof}

\begin{remark}                      
\label{4.17}
It follows from \eqref{4.14a}--\eqref{4.16} that for every $A \in \pi_0 (E
(f, y_0))$ the triple $(C_A, \wtilde{g}_A, \wbar g^\#_A)$ (which
represents $\deg^\#_A (f, y_0)$, cf \eqref{2.6})  contributes just as much
information towards $\deg^\# (f, y_0)$ as the full coincidence data $(C,
\wtilde{g}, \wbar g^\#)$ do. In view of \fullref{1.11} and
\eqref{4.2} this implies
the claim concerning the case $(f_1, f_2) = (f, *)$ in \fullref{1.14} of
the introduction.
\end{remark}

Next recall that -- due to the compactness of $M$ -- our Nielsen numbers
are always finite. Thus according to \fullref{1.14} $N^\# (f, *)$ must
vanish if $b (f, *)$ is infinite. It is rather elementary to show
even more.

\begin{prop}                          
\label{4.18}
If the index of the subgroup $f_* (\pi_1 (M))$ in $\pi_1 (N)$ is
infinite or if $N$ is not compact, then the pair $(f, *)$ is loose;
in other words $f$ is homotopic to a map whose image lies in $N - \{ *\}$.
\end{prop}

\begin{proof} In view of \fullref{1.12} we may assume that $n \ge 2$. Let
$\wtilde{f}$ be a lifting of $f$ to the covering space $\wtilde{N}_f$
(cf \eqref{4.3}) which is noncompact here by assumption. Since $M$ is compact
$\wtilde{f}(M)$ intersects the fiber $p^{- 1} (\{ *\} )$ in only
finitely many points $\wtilde *_1, \dots, \wtilde *_k$. Isotop
$\wtilde{f}$ along disjoint embedded paths $c_i$ in $\wtilde{N}_f$
which start outside of $\wtilde{f} (M)$ and meet the fiber $p^{- 1}
(\{ *\} )$ only at the endpoints $c_i (1) = \wtilde *_i$, $i = 1,
\dots, k$. After this homotopy $\wtilde{f}$ avoids $p^{-1} (\{ *\}
)$, ie $f = p \scirc \wtilde{f}$ maps into $N - \{ *\}$.

(If $N$ is not compact you may also just compose $f$ with a similar
isotopy in $N$ which moves a point $y \in N - f (M)$ to $*$.)
\end{proof}

\begin{proof}[Proof of the root case in \fullref{1.15}]
$f$ and $\wtilde{f}$ allow the same lifting to the common covering
space $\wtilde{N}_f = \wtilde{N}_{\wtilde{f}}$ (compare \eqref{4.3}) of $N$
and $\wtilde{N}$. Thus according to \fullref{4.8} $N^\# (f, *)$ vanishes
if and only if $N^\# (\wtilde{f}, \wtilde * )$ does; otherwise
$$
N^\# (f, *) = [ \pi_1 (N) : f_* (\pi_1 (M))] = d \cdot [\pi_1 (\wtilde
N) : \wtilde{f}_* (\pi_1 (M))] = d \cdot N^\# (\wtilde{f},
\wtilde *)
$$
and similarly for $N (f, *)$. If $d = \infty$, all these Nielsen numbers
vanish.
\end{proof}

\section{Selfcoincidences}                  
\label{sec5}

In this section we prove the results \fullref{1.14} and \fullref{1.15} as far as
selfcoincidences are concerned. Moreover, we relate the selfcoincidence
invariant of a map to its degree.

Let $\underline\Omega^\# (f_1, f_2)$ denote the bordism set of triples
$(C, g, \wbar g^\#)$ as in \eqref{1.6} \eqref{1.7} and \eqref{2.1} ({\it without} a lifting
$\wtilde{g}$). This set is related to $\Omega^\# (f_1, f_2)$ (cf
\eqref{2.1})
via the map $pr_*$ induced by the projection  $pr \co  E (f_1, f_2) \to M$
(cf \eqref{1.6}). In general the resulting coincidence invariant
\begin{equation}                                      
\label{5.1}
\underline\omega^\# (f_1, f_2) := pr_* (\omega^\# (f_1, f_2)) =
[C (f_1, f_2), g, \wbar g^\# ] \in \underline\Omega^\#
(f_1, f_2)
\end{equation}
(compare \eqref{2.2}) is considerably weaker than $\omega^\# (f_1, f_2)$ since
it captures no longer the Nielsen decomposition of the coincidence set,
let alone the other aspects of the lifting $\wtilde{g}$.

Now consider a map $f \co  M \to N$. The selfcoincidence setting is very
special in that $pr$ allows a canonical {\it global} section $s \co  M \to
E (f, f)$ here (defined by $s (x) = (x,$ constant path at $f (x)), \
x \in M$). We obtain induced maps
\begin{equation}                                
\label{5.2}
\xymatrix{
      \underline\Omega^\# (f, f)
      \ar@<0.5ex>[r]^-{s_*}
      &
      \Omega^\# (f, f)
      \ar@<0.5ex>[l]^-{\pr_*}}
\end{equation}
such that $pr_* \scirc s_* =$ identity and clearly
\begin{equation}                              
\label{5.3}
\omega^\# (f, f) = s_* (\underline\omega^\# (f, f))
\end{equation}
(cf \eqref{5.1}; compare also \fullref{2.18}).
Thus only the path component of $E (f, f)$ which contains $s (M)$ can
possibly be strongly essential. Hence if $\omega^\# (f, f) \ne 0$ then
$N^\# (f, f) = 1$. The same argument applies also to the stabilized
coincidence invariant $\wtilde\omega (f, f)$ (which is precisely as
strong as $\omega (f, f) = \pr_* (\wtilde{\omega} (f, f))$, compare
Koschorke \cite{K2,K3}) and $N (f, f)$. In view of \fullref{1.11}(iii)
this proves the claims of \fullref{1.14} as far as they concern the
selfcoincidence case.

As for the proof of \fullref{1.15} just note that $\underline\omega^\#
(\wtilde{f}, \wtilde{f}) = \underline\omega^\# (f, f)$ whenever
$\wtilde{f}$ is a lifting of $f$ into any covering space $\wtilde
N$ of $N$. Indeed, the corresponding coincidence data are essentially
identical, related by the tangent isomorphism of the covering
projection.

Next we show how $\underline\omega^\# (f, f)$ (and hence $\omega^\#
(f, f) = s_* (\underline\omega^\# (f, f)$, cf \eqref{5.3}) can be calculated
once we know $\deg^\# (f) = \omega^\# (f, *)$).

Given $k \in \mathbb Z$ and $[C, \wtilde{g}, \wbar g^\#] \in
\Omega^\# (f, *)$, the homotopy $f|C \sim *$ described by $\wtilde
g$ induces a trivialization of $f^* (TN)|C$ as in \eqref{4.10} and, via
$\wbar g^\#$, of $\nu (C, M)$. As in \fullref{sec4} we may identify
a tubular neighborhood $U$ of $C$ in $M$ with a product $C \times V$
and replace the submanifold $C = C \times \{ 0\} $ by the union
$$C (k) = \bigcup C_i \subset  U \subset M$$
of $|k|$ \lq\lq parallel\rq\rq\ copies $C_i = C \times \{ z_i\} $ where
the points $z_i \in V$, $i= 1, \dots, |k|$, are pairwise distinct. The
obvious linear deformation in $V$ leads to an isomorphism $f^* (TN)|C
\cong f^* (TN)|C_i$; compose it with $\wbar g^\#$ and, if $k$
is negative, with the involution on $C_i \times V \cong \nu (C_i, M)
\cong \nu (C, M)$ determined by a fixed reflection of $V$. We obtain
$$
\wbar g^\# (k) \co  \nu (C (k), M) \cong f^* (TN) | C (k) .
$$
Whenever $n \ge 1$ this construction yields a family of canonical maps
\begin{equation}                           
\label{5.4}
\pinch_{k *} \co  \Omega^\# (f, *) \longrightarrow
\underline\Omega^\# (f, f) , \qquad k \in \mathbb Z,
\end{equation}
which are compatible with homotopies of $f$ and defined by
$$\pinch_{k *} ([ C, \wtilde{g}, \wbar g^\#]) = [C (k) \subset M,
\wbar g^\# (k)].$$

\begin{prop}                      
\label{5.5}
Let $k$ equal the Euler characteristic $\chi (N)$ of $N$ if $N$ is closed
and put $k = 0$ otherwise. Then for every map $f \co  M \to N$ and $* \in N$
we have:
$$
\underline\omega^\# (f, f) = \pinch_{k *} (\omega^\# (f, *)) .
$$
\end{prop}

(The weaker relation $\omega (f, f) = k \deg (f)$ was already proved
and applied in \cite{K2}).

\begin{proof} Identify $V$ with a small neighborhood of $*$ in $N$. After
a homotopy we may assume that
\begin{enumerate}
\item[(i)] $f$ is smooth with regular value $*$, and
\item[(ii)] $C (f, *) = f^{- 1} (\{ *\} )$ has a tubular neighborhood $C
(f, *) \times V  = f^{- 1} (V)$ where $f$ is just the projection to $V$.
\end{enumerate}

Now it is possible to choose a smooth vector field on $N$ with generic
zeroes $z_1, \dots, z_{|k|}$ which all lie in $V$, and without other
zeroes in $f (M)$. Apply the corresponding flow on $N$ to the map $f$
and deform it slightly into a nearby map $f'$. Then the coincidence
locus $C (f', f)$ consists of the inverse images $f^{- 1} (\{ z_i\})$
of the fixed points $z_i$ of the flow, $i = 1, \dots, |k|$. Clearly the
resulting coincidence data represent both $\underline\omega^\# (f, f)$
and $\pinch_{k *} (\omega^\# (f, *))$.
\end{proof}

\begin{example}                            
\label{5.6}
Let $N$ be closed with Euler number $\chi (N) \ne 0$ and  commutative
nontrivial fundamental group. Given $m > n$, consider the projection
$$
f \co   M := N \times S^{m-n} \longrightarrow N .
$$
Then $\underline\omega^\# (f, f) = \pinch_{\chi (N) *} (\omega^\# (f,
*))$ (cf \fullref{5.5}) is nontrivial since already the weaker invariant $\omega
(f, f) = \chi (N) \cdot \omega (f, *)$ (cf \fullref{2.10} and \cite[2.2]{K2})
fails to vanish. This can be detected even by singular homology theory
(which otherwise is often far too crude to capture coincidence phenomena
in higher codimensions). Indeed the composite
$$
\begin{CD}
\Omega^{fr}_{m-n} (M) @>{\mu}>> H_{m-n} (M; \mathbb Z) \
@>{\pi_{2 *}}>> H_{m-n} (S^{m-n}; \mathbb Z) \cong \mathbb Z
\end{CD}
$$
maps $\omega (f, f)$ to $\chi (N)$ (compare the proof of \fullref{5.5}. Note that
$\varphi = f^* (TN) - TM$ is trivial here; $\mu$ and $\pi_2$ denote the
Hurewicz homomorphism and the second projection, resp.).

Thus according to \fullref{1.14}, \eqref{1.9}, and \fullref{3.5}(i) \
$$
N (f, f) = N^\# (f, f) = \MCC (f, f) =  1 \quad\ne\quad \# \pi_0
(E (f, f)) = \# \pi_1 (N)
$$
and $\MC (f, f) = \infty$ whenever $n \ge 1$.
\end{example}

\section{Spherical maps}                       
\label{sec6}

In this section we study in detail the special case $M = S^m$ \. Since the
minimum numbers $\MC$  and $\MCC$ as well as the Nielsen numbers are
(free) homotopy invariants we may -- whenever need be -- assume the maps
$f_1, f_2, \dots$ to have a convenient base point behavior (compare also
\fullref{secA}).

First assume $n \ge 2$ so that we can apply \fullref{2.11} to the case
where $M = S^m$.

Fix basepoints $x_0 \in S^m$ and $y_0, y_1, y_2 \in N$ such that $y_1
\ne y_2$, and choose a local orientation of $N$ at $y_0$ as well as
paths in $N$ joining $y_0$ to $y_1$ and $y_2$. For any two maps $f_i\co  
(S^m, x_0) \to (N, y_i)$, $i = 1, 2,$ these choices allow us to identify
$\Omega^\# (f_1, f_2)$ with the homotopy set $[S^m, S^n \wedge (\Omega
(N)^+) ] \approx \pi_m (S^n \wedge (\Omega (N)^+ ))$ where $\Omega (N)$
denotes the space of loops in $N$ starting and ending in $y_0$ (cf
\fullref{2.11} and its proof). Thus we have a well defined addition
both at the level of maps and of $\omega^\#$--invariants. We will exploit
their compatibilities.

\begin{prop}                        
\label{6.1}
Assume $n \ge 2$. Given $[f_i], [f'_i ] \in \pi_m (N, y_i)$, $i = 1,
2,$ we have
$$
\omega^\# (f_1 + f'_1, f_2 + f'_2 ) = \omega^\# (f_1, f_2) + \omega^\#
(f'_1, f'_2 ) .
$$
In particular,
$$
\omega^\# (f_1, f_2) = \deg^\# (f_1) + \omega^\# (y_1, f_2) .
$$
Furthermore
$$
\deg^\# := \omega^\# (-, y_2) \co  \pi_m (N, y_1) \longrightarrow
\pi_m (S^n \wedge (\Omega N)^+ )
$$
and
$$
\omega^\# (y_1, - ) \co  \pi_m (N, y_2) \longrightarrow \pi_m
(S^n \wedge (\Omega N)^+ )
$$
are group homomorphisms which determine each other via the group
isomorphism $\inv$ defined in (the proof of) \fullref{2.11}.  They measure
also the lack of distributivity of $\omega^\#$, eg
$$
\deg^\#  (f_1) = \omega^\# (f_1, f_2) + \omega^\# (f_1, f'_2) -
\omega^\# (f_1, f_2 + f'_2).
$$
\end{prop}

Thus \eqref{2.17} turns out to be a diagram of homotopy groups and of group
homomorphism when $M = S^m$ and $n \ge 2$.

\begin{proof} The coincidence locus $C (f_1 + f'_1, f_2 + f'_2)$ consists
of the parts $C (f_1, f_2)$ and $C (f'_1, f'_2)$ which lie in disjoint
half-spheres of $S^m$; moreover the coincidence data $\wtilde{g}$ and
$\wbar g$ are compatible with this decomposition. This establishes
the additivity of $\omega^\#$. Obvious homotopies such as $(f_1, f_2)
\sim (f_1 + y_1, y_2 + f_2)$ and $(f_1, y_2) \sim (f_1 + f_1 - f_1,
f_2 + f'_2 - (f_2 + f'_2))$ imply the remaining claims.
\end{proof}

Next let $\MN$ stand for any of the numerical homotopy invariants $N, N^\#,
\MCC$ or $\MC$ for pairs $(f_1, f_2)$ of maps.

\begin{prop}               
\label{6.2}
Let $m, n \ge 1$. Given $[f_i], [f'_i] \in \pi_m (N, y_i)$, $i = 1, 2$,
we have the inequality
$$
\MN (f_1 + f'_1, f_2 + f'_2) \le \MN (f_1, f_2) + \MN (f'_1, f'_2) ;
$$
equality holds if $\MN (f'_1, f'_2) = 0$.

In particular, if $[f_2] \in \pi_m^{(2)} (N, y_2)$ (ie if there is
an element $[\wbar f_2] \in \pi_m (N, y_1)$ such that $(\wbar f_2,
f_2)$ is loose, cf \eqref{1.18}), then
$$
\MN (f_1, f_2) = \MN (f_1 - \wbar f_2, y_2) .
$$
\end{prop}

\begin{proof} The inequality follows from the decomposition
$$C (f_1 + f'_1, f_2 + f'_2) = C (f_1, f_2) \amalg C (f'_1, f'_2).$$
Furthermore note that $\MN (- f'_1, - f'_2) = \MN (f'_1, f'_2)$.
\end{proof}

\begin{cor}                
\label{6.3}
If the homomorphism $i_* \co  \pi_m (N - \{ *\} ) \to \pi_m (N)$ (induced
by the inclusion) is onto, then $\MN (f_1, f_2) = 0$ for all maps $f_1,
f_2 \co  S^m \to N$, ie $(f_1, f_2)$ is loose.

$i_*$ is onto eg if $N \ne S^1$ has an infinite fundamental group,
or if $N$ is not compact, or if $N$ is the product of two manifolds
of strictly positive dimensions, or if $N$ fibers over a manifold $B$
such that $\pi_m (B - \{ *\} ) \to \pi_m (B)$ is onto.
\end{cor}

\begin{proof} After a deformation $f_1$ and $f_2$ have the correct base
point behavior. Then $\MN (f_1, f_2) \le \MN (f_1, y_2) + \MN (y_1, f_2)
= 0$.

If $[f] \in \pi_m (N)$ and $N = N_1 \times N_2$, deform the two
component maps of $f$ until they are constant on opposite half-spheres
$S^m_\pm \subset S^m$; then $f (S^m) \subset N_1 \vee N_2 \subset N -
\{ *\}$. The remainder of our second claim follows from \fullref{4.18}
and the homotopy lifting property of fibrations.
\end{proof}

Propositions \ref{6.1} and \ref{6.2} allow us also in various other situations to
reduce general coincidence questions to the root case (cf \fullref{sec4}) and,
in particular, to a discussion of the degree homomorphism $\deg^\#$.

\begin{proof}[Proof of \fullref{1.20}] We may replace $f_1, f_2$ by
basepoint preserving maps. In view of symmetry results such as \fullref{1.11}(ii)
we may also assume the existence of a loose pair $(\wbar f_2, f_2)$
as in \fullref{6.2} and apply \fullref{1.14} to $f = f_1 - \wbar f_2$. Thus $N^\#
(f_1, f_2) = N^\# (f_1 - \wbar f_2, y_2)$ (cf \fullref{6.2}) equals the
order of $\pi_1 (N)$ when $m \ge 2$ (and of the cokernel of $f_{1 *}
-  f_{2 *}$ when $m = n = 1$, cf \fullref{1.12}) which must be finite
here since $\deg^\# (f_1 - \wbar f_2) = \omega^\# (f_1, f_2) \ne 0$
(cf \fullref{4.18}).

The calculation of $\MCC (f_1, f_2)$ and $\MC (f_1, f_2)$ follows
similarly from \fullref{6.2} and \fullref{1.11}(iii). Analogous conclusions are valid for
$\wtilde{\omega} (f_1, f_2)$ and $N (f_1, f_2)$. Our claims are still
valid for $n = 2$. Indeed, in view of \fullref{6.3}, \fullref{1.26}, and
\fullref{3.5}(iii) we have
to check only the root case when $m > 2$ and $N = S^2$ (or $\mathbb R P
(2)$); after performing a connected sum operation we see that $f^{- 1}
(\{ y_2\} )$ consists of one (or two \lq\lq parallel\rq\rq ) connected
submanifold(s) of $S^m$ (cf \eqref{4.15}).
\end{proof}

How do our minimum numbers behave when $\omega^\# (f_1, f_2)$
vanishes? Let us focus on this question first in the root case.

Given $m \ge 1$, consider the subgroups (cf \fullref{6.1}, \fullref{6.3})
$$i_* (\pi_m (N - \{ *\} )) \subset \ker (\deg^\# ) \qquad \text{of }
\pi_m (N)$$
represented by maps $f \co  S^m \to N$ such that $\MC (f, *) = \MCC (f, *)
= 0$ and that $N^\# (f, *) = 0$, resp. (cf \fullref{4.8}).

\begin{definition}                     
\label{6.4}
$X_m (N) := \ker (\deg^\#) / i_* (\pi_m (N - \{ *\} ))$.

This quotient is always an abelian group. If $n \ge 2$ it measures to
what extend the sequence
\begin{equation}                             
\label{6.5}
\begin{CD}
\pi_m (N - \{ *\} ) @>{i_*}>> \pi_m (N) @>{\deg^\#}>>
\pi_m (S^n \wedge (\Omega N)^+ )
\end{CD}
\end{equation}
fails to be exact. For a description of $X_m (N)$ in terms of pinching
maps see \fullref{7.8} below.
\end{definition}

\begin{thm}               
\label{6.6}
{\rm (a)}\qua Let $Q$ be a smooth connected $q$--manifold and let $p \co  Q \to
N$ be a smooth locally trivial fibration. Assume that the fiber $F =
p^{- 1} (\{ y_0\} )$ is compact (and nonempty) and lies in the
interior of a smoothly embedded $q$--ball $B$ in $Q$. (This holds,
in particular, if $Q$ is a finite covering space over $N)$. Moreover
assume that $q \le 2 n - 2$. Then, given $m \ge 1$, we have: $X_m (N)$
vanishes if and only if $X_m (Q)$ does.

{\rm (b)}\qua $X_m (N) = 0$ in each of the following cases (where $m \ge 1$
is arbitrary unless specified otherwise):
\begin{enumerate}
\item[{\rm (i)}] $m \le 2n - 3$;
\item[{\rm (ii)}] $n \le 2$ (or $m \le 3)$;
\item[{\rm (iii)}] $N$ is not compact;
\item[{\rm (iv)}] $N$ is a sphere $S^n$ or a projective space $\mathbb KP
(n')$, $\mathbb K = \mathbb R, \mathbb C$ or $\mathbb H$, $n, n' \ge 1$;
\item[{\rm (v)}] $N$ is the total space of a Serre fibration with a
section and with strictly positive dimensions of the fiber and base space;
\item[{\rm (vi)}] $N$ fibers over a manifold which has an infinite
fundamental group.
\end{enumerate}
\end{thm}

\begin{proof} \fullref{1.12} establishes these claims whenever $m = 1$ or
$n = 1$. Thus we may assume $m, n \ge 2$ in the remainder of this proof.

Given a smooth fibration $p \co  (Q^q, \wtilde{y}_0) \to (N^n, y_0)$,
consider the diagram
\begin{equation}                                
\label{6.7}
  \xymatrix{
    \pi_m(Q,\wtilde{y}_0)
    \ar[r]^-{p_*} \ar[d]^-{\operatorname{deg}_Q^\#}  &
    \pi_m(N,y_0)
    \ar@<0.5ex>[r]^-{\partial}  \ar[d]^-{\operatorname{deg}_N^\#} &
    \pi_{m-1}(F,\wtilde{y}_0)
    \ar@<0.5ex>@{-->}[l]^-{j_*}  \\
    \pi_m(S^q\wedge(\Omega Q)^+)
    \ar@<-0.5ex>@{-->}[r]_-{\beta}  &
    \pi_m(S^n\wedge(\Omega N)^+)
    \ar@<-0.5ex>[l]_-{\alpha}
  }
  \end{equation}
which involves the homotopy sequence of $p$. The homomorphism $\alpha$
is defined as follows. Interpret an element $c \in \pi_m (S^n \wedge
(\Omega N)^+)$ -- via Pontryagin--Thom -- as a bordism class of a framed
submanifold $C \subset \mathbb R^m$, together with a map $\wtilde{g}
\co  C \to \Omega (N, y_0)$. The corresponding evaluation map $C \times I
\to N$ lifts to a homotopy $\wtilde{G}$ in $Q$ from the constant map
at $\wtilde{y}_0$ to a map $\wtilde{G}_1 \co  C \to F$ which we may
assume to be smooth, with regular value $\wtilde{y}_0$. Endow $C' :=
\wtilde{G}^{- 1}_1 (\{ \wtilde{y}_0\} )$ with the map $\wtilde{g}'
\co  C' \to \Omega (Q, \wtilde{y}_0)$ which corresponds to $\wtilde{G} |
C' \times I$. Moreover compose the natural framing of $C'$ in $C$ (given
by tangent map of $\wtilde{G}_1$) with the automorphism of $C' \times
T_{\wtilde{y}_0} F$ which is determined by the homotopy $\wtilde
G | C' \times I$ and the tangent bundle along the fibers of $p$ (cf
\cite[3.1]{K3}). The resulting bordism class $[C' \subset \mathbb R^m,
\wtilde{g}']$ defines $\alpha (c)$. We have
\begin{equation}                          
\label{6.8}
\deg^\#_Q  = \alpha \scirc \deg^\#_N \scirc p_*
\end{equation}
since $\deg^\#_N \scirc p_*$ and $\deg^\#_Q$ correspond to taking the
inverse image of a fiber and of a point in $Q$, resp.

If there exists a basepoint preserving homotopy $\wtilde{J} \co F \times
I \to Q$ from the constant map at $\wtilde{y}_0$ to the fiber inclusion
then $p \scirc \wtilde{J}$ corresponds to a map $j \co  F \to \Omega N$;
this induces a splitting of the top line in diagram \eqref{6.7} since $\partial
\scirc j_* = \id$. If in addition $q \le 2n - 2$, ie $\dim (F \times
I) < \dim N$, then the image of $j_*$ lies already in $i_* (\pi_m (N -
\{ *\}, y_0))$ for $* \ne y_0$.

Finally assume the full hypothesis in \fullref{6.6}(a)). Given a framed
$q$--codimensional submanifold  $C' \subset \mathbb R^m$ together with a
map $\wtilde{g}' \co  C' \to \Omega (Q, \wtilde{y}_0)$, twist the framing
via the automorphism of $C' \times \mathbb R^q$ determined by $\wtilde
g'$ and $TQ$, use the twisted framing to identify a tubular neighbourhood
of $C'$ in $\mathbb R^m$ with $C' \times B$, reframe the submanifold
$$C := C' \times F \subset C' \times B \subset \mathbb R^m$$
using $p, \wtilde{g}'$ and a contraction $\wtilde{J} \co  F \times I
\to B \subset Q$, and equip $C$ with the paths in $N$ which concatenate
$p \scirc \wtilde{g}'$ with the adjoint $j$ of $p \scirc \wtilde
J$. The resulting bordism class $[C, \wtilde{g}]$ corresponds to $\beta
([C', \wtilde{g}'])$.

This definition of $\beta$  mimics the transition from $\deg^\#_Q$  to
$\deg^\#_N$ where the inverse image of a point $\wtilde * \in Q$
is replaced by the inverse image of a whole fiber $p^{- 1} (\{ *\}
)$ containing $\wtilde *$ (where $* \in N - \{ y_0\}$ is close to
$y_0$). The reframings in the construction of $\alpha$ and $\beta$
are motivated by our framing convention in the definition of $\deg^\#$
(cf \eqref{2.15}). We obtain
\begin{equation}                            
\label{6.9}
\beta \scirc \deg^\#_Q = \deg^\#_N \scirc p_*
\end{equation}
and $\alpha \scirc \beta = \id$. Hence $\beta$ is injective.

Moreover the following three statements are equivalent for every $[f] \in
\pi_m (Q)$: $f$ can be deformed into the complement of (i) $\{\wtilde
*\}$, (ii) $\mathring{B}$, (iii) $F$. This completes the proof
of claim (a) in \fullref{6.6}.

Clearly $X_m (N) = 0$ whenever $i_*$ is onto (eg in case (v) of claim
(b); cf \fullref{6.3} and its proof) or when $m < n$. Moreover $X_m$
vanishes for all spheres (since $\deg^\#$ is split injective here,
cf \fullref{2.18}) and hence also for all projective spaces (by claim (a)). In
particular, $X_m (N) = 0$ for $N = S^2$ or $\mathbb R P^2$ and also
for the remaining surfaces (which are noncompact or have infinite
fundamental groups). This settles case (ii) of \fullref{6.6}(b). Case (i) is a
(weak) consequence of \cite[Theorem~1.10]{K3}.
\end{proof}

\begin{remark}                      
\label{6.10}
In the previous discussion we have been dealing with the very special case
where $M$ is a sphere and $f_2$ is constant. Here no problem arises in the
dimension setting $m = n = 2$ which is so critical in classical fixed
point theory: see Jiang \cite{Ji2} for examples where $N (f, \id) = N^\# (f,
\id) = 0 < \MCC (f, \id) \le  \MC (f, \id)$.
\end{remark}

\begin{example}[Stiefel manifolds and Grassmannians]  
\label{6.11}

(a)\qua  Given integers $1 \le k < r$, let $p_{r, k} \co  V_{r, k} \to G_{r,
k}$ be the fiber projection which maps each orthonormal $k$--frame in
$\mathbb R^r$ to the $k$--plane it spans. Assume $2k \le r$. (This is
no restriction as far as Grassmannians are concerned since $G_{r, k}
\cong G_{r, r-k})$. Then the fiber $O (k) = V_{k, k}$ of $p_{r, k}$
can be deformed into a ball via the isotopy
$$
V_{k,k} \times \bigl[0, \tfrac{\pi}{2} - \varepsilon\bigr] \longrightarrow
V_{r,k},
$$
$(\{v_i\}, t) \to \{ \cos (t) v_i + \sin (t) e_{k +i}\}$, where $1 \le
i \le k$ and $e_j$ denotes the $j^{\text{th}}$ standard unit vector in
$\mathbb R^r$. Thus, given $m \ge 1$, we conclude from \fullref{6.6}(a)) that $X_m
(G_{r, k}) = 0$ if and only if $X_m (V_{r, k}) = 0$.

(b)\qua In turn, if the sphere $S^{r -1}$ allows $k - 1$ linearly independent
vector fields then the canonical fibration $V_{r, k} \to S^{r -1}$
has a section and $X_m (V_{r, k}) = 0$ by \fullref{6.6}(b)(v).

We conclude for instance that $X_m (G_{r, k}) = 0$ for all $m \ge 1$ and
$1 \le k < r$ if $r =  2, 4, 8$, or $16$ (compare eg
Atiyah--Bott--Shapiro \cite[page~38]{ABS}).

The same statements hold for the Grassmann manifold $\wtilde{G}_{r, k}$
of oriented $k$--planes in $\mathbb R^r$ (with the added advantage that no
map into $\wtilde{G}_{r, k}$ is coincidence producing).
\end{example}

\begin{remark}[and proof of \fullref{1.23}]
\label{6.12}
The homogeneity techniques
supplied by Propositions \ref{6.1} and \ref{6.2} allow us (in analogy to the proof
of \fullref{1.20}) to extend the condition \lq\lq $X_m (N) = 0$\rq\rq\ to yield
the following equivalent condition: \lq\lq For every pair of maps $f_1,
f_2 \co  S^m \to N$ which are not both coincidence producing we have: $(f_1,
f_2)$ is loose if and only if $\omega^\# (f_1, f_2) = 0$\rq\rq.

Clearly this implies \fullref{1.23} of the introduction.
\end{remark}

\begin{proof}[Proof of the statements in \fullref{1.24}] In view of
\fullref{1.12} we
may assume that $n \ge 2$. In the case $N = S^n$ $\deg^\#$ is split
injective (cf \fullref{2.18}) and $X_m (N) = 0$ (cf \fullref{6.6}). Moreover $(f_1, f_2)$
is loose if and only if $f_1 \sim a f_2$ (cf Dold--Gon\c{c}alves \cite[1.10]{DG}
or Koschorke \cite[Section~8]{K3}). Thus it follows from \fullref{1.20}
and \fullref{1.23} that $N^\# (f_1,
f_2) = \MCC (f_1, f_2)$ equals $0$ or $\# \pi_0 (E (f_1, f_2))$
according as $\omega^\# (f_1, f_2) = \deg^\# (f_1 - a \scirc f_2)$
(cf \fullref{6.1}) vanishes or not. The restriction $n \ne 2$ in
\fullref{1.20} can be
avoided here by a connected sum argument applied to the (generically)
framed manifold $(f_1 - a \scirc f_2)^{- 1} \{ *\}$ (compare
\fullref{2.18}).

If $\MC (f_1, f_2)$ is finite and $m > n \ge 2$ then $\deg^\# (f_1 -
a \scirc f_2) = e (\beta)$ for some $\beta \in \pi_{m -1} (S^{n -1})$
(cf \fullref{3.5}) and hence in view of \fullref{2.18} $[f_1 - a \scirc f_2] = \coll_*
([ f_1 - a \scirc f_2]) = \pm E (\beta)$.
\end{proof}

Before we discuss the minimum number $\MC (f_1, f_2)$ for maps into
a general target manifold $N$ we introduce the composite homomorphism
(for $m > 2 \le n)$
\begin{equation}                    
\label{6.13}
\begin{CD}
\eta \co  \bigoplus_{A \in \pi_1 (N)} \pi_{m -1} (S^{n -1})
\longrightarrow \pi_{m -1} \Bigl(\bigvee_A S^{n -1} \Bigr) @>{ v|_* }>> \pi_{m
-1} (N - \{ *\} )
\end{CD}
\end{equation}
as follows. Choose a homeomorphism $u$ from the unit ball $B^n \subset
\mathbb R^n$ onto a small ball $B (*)$ in $N$ centered at $*$ such
that $y_0 = u (x_0) \in \partial B (*)$ (where $x_0 \in S^{n -1}$ is a
basepoint); lift $u$ in the universal cover $\wtilde{N}$ to all the
various levels (which are parametrized by $A \in \pi_1 (N, y_0) \approx
\pi_0 (E (f_1, f_2))$); join the lifted basepoints $\{\wtilde{u}_A
(x_0)\}$ by appropriate paths to $\wtilde{y}_0$ and project back to
$N$. We obtain the map
\begin{equation}                     
\label{6.14}
v \co  \Bigl( \bigvee_A B^n, \bigvee_A S^{n -1} \Bigr) \longrightarrow
(N, N - \{ *\} ) .
\end{equation}
Compose the induced homomorphism $v|_*$ with the natural inclusion to
define $\eta$ in \eqref{6.13}.

On the other hand consider the maps
$$
\xymatrix{
     \bigvee\limits_{A \in \pi_1 (N)}^{~} S^n
      \ar@<0.5ex>[r]^-{\In}
      &
      S^n \wedge (\Omega N^+) =  \text{Thom} (\Omega N \times
      \mathbb R^n)
      \ar@<0.5ex>[l]^-{\ret}
  }
$$
corresponding (at the level of Thom spaces of trivial $n$--plane bundles)
to the obvious componentwise inclusion and retraction maps between $\pi_1
(N) = \pi_0 (\Omega N)$ and $\Omega N$. We obtain the diagram (for $m
\ge 2$)
\begin{equation}                                 
\label{6.15}
\xymatrix{
     \bigoplus\limits_{A \in \pi_1 (N)} \pi_m (S^n) \subset \pi_m
     \Bigl( \bigvee\limits_A S^n\Bigr)
     \ar@<0.5ex>[r]^-{\In_*}
     &
     \pi_m (S^n \wedge (\Omega N)^+)
     \ar@<0.5ex>[l]^-{\ret_*}
   }
\end{equation}
where $\ret_* \scirc \In_* = $ identity. Moreover, $\In_*$ and the
Freudenthal suspension compose to yield the homomorphism
\begin{equation}                          
\label{6.16}
\pm e = \In_*| \scirc \oplus E \co  \bigoplus_{A \in \pi_1 (N)}
\pi_{m -1} (S^{n -1}) \longrightarrow \pi_m (S^n \wedge (\Omega N)^+ )
\end{equation}
constructed at the beginning of Sections \ref{sec3} and \ref{sec6}.

\begin{thm}                            
\label{6.17}
Assume $m \ge 3$. Given any pair of maps $f_1, f_2 \co  S^m \to N$ such
that $f_2$ is not coincidence producing, we have:

If $\MC (f_1, f_2) < \infty$ then $\omega^\# (f_1, f_2) \in e (\ker
(\eta))$. In turn, if $X_m (N) = 0$ and $\omega^\# (f_1, f_2) \in e
(\ker (\eta))$ then $\MC (f_1, f_2)$ is finite.
\end{thm}

\begin{proof} In view of Propositions~\ref{6.1} and~\ref{6.2} we need to consider only pairs
of the form $(f, *)$. Interpret $f$ as a map from $I^m$ to $N$ which
maps the boundary $\partial I^m$ to the basepoint $y_0 \ne *$. If
$\MC (f, *) < \infty$ then after a deformation $f^{- 1} (\{ *\} )$
consist of finitely many points $x_j \in \mathring{I}{^m}, \
j = 1, \dots, k$. Furthermore there are small balls $B_{x_j} \subset
\mathring{I}{^m}$ centered at $x_j$ such that $f$ maps $\partial
B_{x_j} \cong S^{m -1}$ to the boundary sphere of the ball $B (*)$
around $*$ in $N$. After isotoping $x_j$ into suitable positions in
$I^m = I^{m -1} \times I$ we may even assume that
$$
\begin{aligned}
f & \bigl(I^{m -1}\times\bigl[0,\tfrac12\bigr]\bigr) \subset v \Bigl(\bigvee_A
B^n\Bigr) \quad \text{(cf \eqref{6.14})} \quad \text{and}  \\
f & \bigl(I^{m -1} \times \bigl[\tfrac12, 1\bigr]\bigr) \subset N - \mathring{B}
(*).
\end{aligned}
$$
Thus the index maps of the coincidence points of $(f, *)$ (cf
\eqref{3.4})
determine an element $\alpha \in \ker \eta$ (since $f | I^{m-1} \times
\bigl\{\frac12\bigl\} $ is nulhomotopic in $N - \{ *\} ) $ such that $e
(\alpha) = \omega^\# (f, *)$.

On the other hand, every $\alpha \in \ker (\eta)$ and a corresponding
nulhomotopy in $N - \{ *\}$ yields a map $f' \co  S^m \to N$ with $\MC
(f', *) < \infty$ and $\omega^\# (f', *) = e (\alpha)$. If $\omega^\#
(f, *) = e (\alpha)$ and $X_m (N)$ vanishes then so do $\deg^\# (f -
f')$ and hence $\MC (f - f', *)$.

Since $\MC (f, *) \le \MC (f', *) + \MC (f - f', *)$ we conclude that
this minimum number is finite.
\end{proof}

We illustrate our criterion by a sample application.

\begin{cor}                               
\label{6.18}
Let $N = S^n / G$ be a spherical space form (compare \fullref{1.25}) and assume
$m, n \ge 2$. Then we have for all $[f] \in \pi_m (N)$
$$
\MC (f, *) < \infty \Longleftrightarrow [\wtilde{f}]
\in
\begin{cases}
E (\ker h) &\text{if } \# G \ge 3 ; \\
E (\pi_{m -1} (S^{n -1})) &\text{if } \# G \le 2.
\end{cases}$$
Here $[\wtilde{f}] \in \pi_m (S^n)$ is obtained by lifting $[f]$, $E$
denotes the Freudenthal suspension and
$$
h := \bigoplus^\infty_{j = 0} h_j \co  \pi_{m -1} (S^{n -1})
\longrightarrow \pi_{m -1} (S^{2n -3}) \oplus \pi_{m -1} (S^{3n -5})
\oplus \dots
$$
is the (total) Hopf--Hilton homomorphism (cf Whitehead \cite[XI, 8.5]{W}.
\end{cor}

Recall that $E (\pi_{m -2} (S^{n -2}))$ lies  in the kernel of $h$
whenever $n \ge 3$ (and coincides with it if also $m \le 3 n - 6$; see,
for example, \cite[XII, 2.3]{W} or Koschorke--Sanderson \cite{KS}).

\begin{proof} If $n$ is even, then $\chi (S^n) = 2 \ge \# G$. In
particular, this implies our claim when $m = n = 2$.

Thus (in view of \fullref{3.5}(iii))  we may assume that $m, n \ge 3$. Then $\eta$
(cf \eqref{6.13}) lifts to
$$\wtilde\eta \co (\pi_{m -1} (S^{n -1}))^k \subset \pi_{m -1}
  \Bigl(\bigvee\limits^k\!\! S^{n-1}\Bigr)
\stackrel{\wtilde{v} |_*}{\longrightarrow}
  \pi_{m -1} (S^n {-}\{ \wtilde *_1\!,\!.,\!\wtilde *_k\}) \cong \pi_{m -1}
  \Bigl(\bigvee\limits^{k-1}\!\!S^{n -1}\Bigr)$$
where $k := \# \pi_1 (N) = \# G$ and the points $\wtilde *_1, \dots,
\wtilde *_k \in S^n$ are projected to $* \in N$. After a suitable
homotopy $v | \bigvee^{k -1} S^{n -1}$ (cf \eqref{6.14}) lifts to an embedding
$$
\wtilde{v} \Vert \co \bigvee^{k -1} S^{n -1} \subset \mathbb R^n
- \{ \wtilde *_1, \dots, \wtilde *_{k -1} \} \approx S^n - \{
\wtilde *_1, \dots, \wtilde *_k\}
$$
which is known to be a homotopy equivalence. We use this to identify
the target group of $\wtilde\eta$ with $\pi_{m -1} \Bigl(\bigvee^{k -1}
S^{n -1}\Bigr)$ and thus make it accessible to Hilton's computing techniques
\cite{H}.

Let $\iota_j \in \pi_{n -1} \Bigl(\bigvee S^{n -1}\Bigr)$ be represented by the
inclusion of the $j$th sphere in the wedge. Then
\begin{equation}                                      
\label{6.19}
\begin{aligned}
\wtilde{v} |_* (\iota_j) &= \iota_j, &j = 1, \dots, k - 1, \text{ and} \\
\wtilde{v} |_* (\iota_k) &= - \iota_1 -  \cdots - \iota_{k-1};&
\end{aligned}
\end{equation}
indeed, we can retract the boundary sphere of a ball around the point
$\wtilde *_k = \infty \in \mathbb R^n \cup \{ \infty\}$ into the
spheres around the remaining points $\wtilde *_1, \dots, \wtilde
*_{k -1}$, reversing orientations in the process. Let $r$ stand for a
suitable reflection (ie selfmap of degree $- 1$) on any sphere.

Now consider an element $\alpha = (\alpha_1, \dots, \alpha_k)$ in the
domain of $\wtilde\eta$. When $k = 2$ then $\eta (\alpha) = 0$ if
and only if $\alpha_2 = -r_* (\alpha_1)$; in this case $E (\alpha_1) =
E (\alpha_2)$. When $k = 3$ then $\eta (\alpha)$ vanishes if and only if
\begin{equation}                               
\label{6.20}
(\iota_1 + \iota_2) \scirc \alpha_3 = \iota_1 \scirc \alpha_3 +
\iota_2 \scirc
 \alpha_3 + \sum_{j \ge 0} w_j (\iota_1, \iota_2) \scirc h_j (\alpha_3)
 \end{equation}
 (cf Whitehead \cite[XI, 8.5]{W}) equals $- (- \iota_1) \scirc \alpha_1 - (-
 \iota_2) \scirc \alpha_2$, or, equivalently, $\alpha_1 = \alpha_2 =
 \alpha_3 \in \ker h$. (In order to see this, project to each of the
 wedge factors $S^{n -1}$ and permute their roles; note also  that $\ker
 h \subset \ker (r_* + \id)$: just substitute $\iota_1$ and $\iota_2$,
 resp., by $r$ and $\id$, resp., in \eqref{6.20}). Collapsing all but
 three spheres in the wedge $\bigvee S^{n -1}$ allows us to extend our
 calculation of $\ker \eta = \ker \wtilde\eta$ also to the case where
 $k > 3$.

Next observe that $X_m (N) = 0$ (cf \fullref{6.6}(a) and (b)(iv)). Thus
according to \fullref{6.17} $\MC (f, *)$ is finite if and only if $\deg^\#
(f) \in e (\ker \eta)$ (cf \eqref{6.16}) or, equivalently, $\deg^\# (\wtilde
f)$ lies in the image, under $\In_{\ell *} \scirc E$ for some suitable
loop $\ell \in \Omega (N)$ (cf \eqref{2.17}), of $\pi_{m -1} (S^{n -1})$ and
of $\ker h$, resp. (see \fullref{4.17}). But this yields the indicated
condition concerning the homotopy class  $[\wtilde{f}] = \coll_*
([\wtilde{f}])$ (cf \fullref{2.18}).
\end{proof}

\begin{example}                       
\label{6.21}
Let $N = S^3 / G$ be a 3--dimensional spherical space form and $m =
4$. Then $\begin{CD} h \co  \pi_3 (S^2) @>{\cong}>> \mathbb Z\end{CD}$
and $E$ maps this group onto $\pi_4 (S^3) \cong \pi_4 (N) \cong \mathbb
Z_2$. Hence according to \fullref{1.14}, \fullref{2.18}, \fullref{4.8},
and \fullref{6.18} we have for every
map $f \co  S^4 \to N$
$$\MC (f, *) =
\begin{cases}
\infty & \text{if } [f] \ne 0 \text{ and } \# G \ge 3; \\
\# G & \text{if } [f] \ne 0 \text{ and }  \# G \le 2; \\
0 & \text{if } [f] = 0.
\end{cases}$$
In particular, if $\# G \ge 3$ and $[f] \ne 0$ then
$$\MC (f, *) = \infty \quad \text{but} \quad
  \MC (\wtilde{f}, \wtilde *) = 1$$
where $\wtilde{f} \co  S^4 \to S^3$ is a lifting of $f$.
\end{example}

In order to gain a better geometric understanding of the dependence on
$\# G$ which we encountered in \fullref{6.18} and in the last example let
us recall (and generalize) a classical notion (cf Wyler \cite[page~29]{Wy}).

\begin{definition}            
\label{6.22}
Given a map $f \co  M \to N$ we call a point $y \in N$ {\it injective}
(or {\it almost injective}, resp.) if its inverse image $f^{-1} (\{ y\}
)$ consists of a single point (or is finite, resp.).
\end{definition}

Thus $\MC (f, *) < \infty$ if and only if $f$ is homotopic to a map
which possesses an almost injective point. Furthermore, given a lifting
$\wtilde{f}$ of $f$ to a $k$--sheeted covering space $\wtilde{N}$
over the manifold $N$, $\MC (f, *) < \infty$ if and only if $\wtilde{f}$
is homotopic to a map which has at least $k$ different almost injective
points (which we may isotop into the fiber $\{ \wtilde *_1, \dots,
\wtilde *_k\}$ over $* \in N$).

Clearly if a map between spheres is an (unreduced) suspension then it
has a least two injective points. Suspending again, we obtain a map which
admits a whole circle of injective points. In \fullref{6.18} and
\fullref{6.21} Hopf--Hilton invariants describe the precise borderline between
single and double suspensions which characterizes spherical maps having
at least as many almost injective points as there are elements in $G$.

\begin{proof}[Proof of the statements in \fullref{1.25}] Since $n$ is odd
$N$ admits a nowhere vanishing vector field. Hence for every map $f$
into $N$ the pair $(f, f)$ is loose. In particular if $f_1 \not\sim f_2$
and $m > 1$ then $\omega^\# (f_1, f_2) = \deg^\# (f_1 - f_2)$ (cf
\fullref{6.1})
is nontrivial since $X_m (N) = 0$ (cf \fullref{6.6}(a) and (b)) and we see from
a lifting argument that $f_1 - f_2$ does not map into $N - \{ *\}$. The
statements in \fullref{1.25} follows from \fullref{1.20} and
\fullref{6.18} (see also \fullref{1.12}).  \end{proof}

\begin{remark}                             
\label{6.23}
The approach of the preceding proof still works when $n$ is even  {\it
provided} $f_2$ occurs in a loose pair $(\wbar f_2, f_2)$ (so that
we may replace $(f_1, f_2)$ by $(f_1 - \wbar f_2, *)$, cf
Propositions~\ref{6.1} and \ref{6.2}). But this proviso can be very
restrictive. For instance if $N$ is
an even-dimensional real projective space and $m \ge 2$, then it implies
that both liftings $\wtilde{f}_2$ and $a \scirc \wtilde{f}_2$ of $f_2$
in $S^n$ must be homotopic to $a \scirc \overset{\simeq}{f}_2$ where $a$
denotes the antipodal map (cf Dold--Gon\c{c}alves \cite[1.10]{DG}).
Thus $\wtilde{f}_2
\sim a \scirc \wtilde{f}_2$ and, if $[\wtilde{f}_2] \in E (\pi_{m -1}
(S^{n -1}))$, then $2 [f_2] = 0$. This restriction is satisfied, for
example,
by only two elements in $\pi_{17} (\mathbb R P (10)) \cong \pi_{17}
(S^{10}) = E (\pi_{16} (S^9) \cong \mathbb Z_{240}$.
\end{remark}

Next consider again an arbitrary manifold $N$ as well as a pair of maps
$f_1, f_2 \co  S^m \to N$, $m \ge 2$, with only finitely many coincidence
points. Adding up their index maps (cf \eqref{3.4}) while keeping track of
the Nielsen decomposition we obtain the \lq\lq total index\rq\rq\
$$
\ind (f_1, f_2) = \bigl\{ \ind_A (f_1, f_2)\bigr\} \in \bigoplus_{A
\in \pi_1 (N)} \pi_{m -1} (S^{n -1})
$$
which is well defined (once the choices of an orientation of $N$ at $y_0$
and of paths joining $y_0$ to $f_1 (x_0)$ and $f_2 (x_0)$ have been
fixed as at the beginning of this section).

\begin{thm}                    
\label{6.24}
The suspended index
$$\bigoplus E (\ind (f_1, f_2)) \in \bigoplus_{A \in \pi_1 (N)} \pi_m (S^n)$$
is determined by $\omega^\# (f_1, f_2)$ and hence depends only on the
basepoint preserving  homotopy classes of $f_1$ and $f_2$.

Furthermore, if the suspension $E \co  \pi_{m -1} (S^{n -1}) \to \pi_m
(S^n)$ is injective  then
$$\MC (f_1, f_2) = N^\# (f_1, f_2) =
  \# \{ A \in \pi_1 (N) | \ind_A (f_1, f_2) \ne 0\}.$$
\end{thm}

\begin{proof} When we apply the monomorphism $\In_*$ (cf \eqref{6.15}) to the
suspended index we obtain $e (\ind (f_1, f_2)) = \omega^\# (f_1, f_2)$
(cf \eqref{6.16} and the discussion preceding \fullref{3.5}).

If $n > 2$ we may replace any Nielsen class by (at most) one single point
$x$ (cf \fullref{3.5}(ii)); if it is nonessential we may remove it altogether
provided $E$ is injective and hence the index map $q$ at $x$ (cf
\fullref{3.4})
is nulhomotopic.  In view of \fullref{3.5}(iii), \fullref{6.2},
\fullref{6.3}, and (the last statement
in) \fullref{2.18} it remains only to consider the case where $m = 2$ and $N =
\mathbb R P (2)$. But here our second claim follows from Jezierski
\cite[Theorem~4.0]{Je}. (I am grateful to S Bogaty\u\i\ and E Kudryavtseva
for bringing this reference to my attention).
\end{proof}

\section{Hopf--Ganea invariants and the degree $\deg^\#$}
\label{sec7}

In this section we continue our discussion of the case $M = S^m$,
$m \ge 2$. We give a purely homotopy theoretical description of our
(geometric) homomorphism $\deg^\#$ and of the group $X_m (N)$. $\deg^\#$
turns out to consist of two components, one of them being a Hopf--Ganea
invariant $H_{\mathcal C}$. We show that $H_{\mathcal C}$ is a finiteness
obstruction for the minimum number $\MC$ (and actually the only one in
a dimension range depending on the connectivity of $N$).

Assume that $k := \# \pi_1 (N)$ is finite (otherwise $\deg^\# \equiv 0$,
cf \fullref{4.18}). Let $p \co  \wtilde{N} \to N$ be the universal covering of
$N$ and let $\bigvee B^n$ denote a wedge of standard unit $n$--balls,
pinched together at base points of the boundary spheres. Given points
$y_0 \ne * \in N$ and $\wtilde{y}_0 \in p^{- 1} (\{ y_0\} )$ as well
as an orientation of $\wtilde{N}$, we can identify $\bigvee^k B^n$ with
a union of (compact) balls around $\wtilde *_1, \dots, \wtilde *_k
\in p^{- 1} (\{ *\} )$ in $\wtilde{N}$ which intersect at $\wtilde
y_0$. We obtain a pinching map
\begin{equation}                                    
\label{7.1}
\pinch\co \wtilde{N} = (\wtilde{N} - \cup^k
\mathring{B}{^n} ) \cup_{\att} \bigvee^k B^n \longrightarrow
\wtilde{N} / \bigvee^k S^{n -1} \cong \bigvee^k S^n \vee \wtilde{N}
\end{equation}
which collapses the boundary spheres $\bigvee S^{n -1}$ to a point.

On the other hand there is the projection
\begin{equation}                                          
\label{7.2}
\proj_2\co  \bigvee^k S^n \vee \wtilde{N} \longrightarrow
\wtilde{N}
\end{equation}
which collapses the wedge $\bigvee ^k S^n$. Its mapping fiber (cf
Whitehead \cite[page~43]{W} but  {\it without} a retopologization as in
\cite[page~20, line~7]{W}).
\begin{equation}                                        
\label{7.3}
F = \Bigl\{ (y, \theta) \in \Bigl(\bigvee^k S^n \vee \wtilde{N}\Bigr) \times P
(\wtilde{N}) ~|~ \theta (0) = \proj_2 (y), \theta (1) = \wtilde
y_0\Bigr\}
\end{equation}
(compare \fullref{1.27}) contains the contractible subspace
$$
P = \{ (y, \theta) \in F ~|~ y \in \wtilde{N}\} \cong \{
\theta \in P (\wtilde{N}) ~|~ \theta (1) = \wtilde{y}_0 \}
$$
whose complement $(S^n - \{ \wtilde{y}_0\} ) \times (\amalg_k \Omega
(\wtilde{N}, \wtilde{y}_0))$ admits a homotopy equivalence with $(S^n
- \{ \wtilde{y}_0\} ) \times \Omega (N, y_0)$ which is compatible with
the labelling of the center points $\wtilde *_1, \dots, \wtilde
*_k \in p^{- 1} (\{ *\} )$ of the balls in $\bigvee^k B^n$.

\begin{lem}[cf Cornea {{\cite[page~2769]{C2}}}]   
\label{7.4}
The quotient map
$$\quot \co F \longrightarrow F / P \cong S^n \wedge (\Omega(N, y_0)^+)$$
is a homotopy equivalence.
\end{lem}

\begin{proof} Let $\theta_0$ denote the constant loop at $\wtilde
y_0$. We describe a map
$$
\Bigl( \Bigl(\bigvee^k S^n\Bigr) \times \Omega (\wtilde{N}, \wtilde{y}_0), \{
\wtilde{y}_0\} \times \Omega (\wtilde{N}, \wtilde{y}_0)\Bigr) \
\longrightarrow (F, \{ (\wtilde{y}_0, \theta_0)\} )
$$
which induces the required homotopy inverse. In each sphere $S^n$
contract all points $y$  with  $\dist (y, \wtilde{y}_0) = 1$ into the
wedge point $\wtilde{y}_0$ and use this distance also as a parameter
for how far to deform loops $\theta \in \Omega (\wtilde{N}, \wtilde
y_0)$ along themselves. Thus at $y = \wtilde{y}_0$ each loop is fully
retracted to $\theta_0$ while $\theta$ is left unchanged if $\dist (y,
\wtilde{y}_0) \ge 1$.
\end{proof}

The maps in \eqref{7.1}, \eqref{7.2}, and \fullref{7.4} yield the diagram
\begin{equation}
\label{7.5}
  \xymatrix{
     \pi_m(N,y_0)
     \ar@{-->}[d]^-{\deg^\#}
     &
     \pi_m(\wtilde{N},\wtilde{y}_0)
     \ar@{-->}[l]_-{p_*}^-{\cong}
     \ar@{-->}[ld]^-{H}
     \ar[d]^-{\operatorname{pinch}_*}
     \\
     \!\!\!\!0 \to \pi_m(S^n{\wedge}(\Omega(N,y_0)^+))
     \ar@<0.5ex>[r]^-{i_{1*}}
     &
     \pi_m\Bigl(\bigvee\limits^k_{~} S^n{\vee}\wtilde{N},\wtilde{y}_0\Bigl)
     \ar@<0.5ex>@{-->}[l]^-{\operatorname{pr}_1}
     \ar@<0.5ex>[r]^-{\operatorname{proj}_{2*}}
     &
     \pi_m(\wtilde{N},\wtilde{y}_0)
     \ar@<0.5ex>@{-->}[l]^-{\operatorname{incl}_{2*}}
     \to
     0 }
\end{equation}
where the horizontal short exact homotopy sequence (of $\proj_2$,
turned into a fiber map) splits canonically via the inclusion $\incl_2\co  
\wtilde{N} \subset \bigvee S^n \vee  \wtilde{N}$. As we will see below
(cf \eqref{7.11}) the composite $H$ of the resulting projection $\pr_1$ with
$\pinch_*$ is an enriched Hopf--Ganea invariant homomorphism (relative to
the attaching map $\att$ in \eqref{7.1}). But let us first compare it
to $\deg^\# = \omega^\# (~,*)$ (cf \fullref{6.1}).

\begin{thm}                      
\label{7.6}
$H := \pr_1 \scirc \pinch_*$ coincides with $\deg^\# \scirc p_*$ up to
an involution of the target group $\pi_m (S^n \wedge \Omega (N, y_0)^+)$.
\end{thm}

\begin{proof} Given $\wtilde{f} \co  (I^m, \partial I^m) \longrightarrow
(\wtilde{N}, \wtilde{y}_0)$, we define
$$u', u'', u \co  (I^m, \partial I^m) \longrightarrow \Bigl(\bigvee^k S^n \vee
  \wtilde{N}, \wtilde{y}_0\Bigr)$$
by $u' := \pinch \scirc \wtilde{f}$, $u'' = \incl_2 \scirc \proj_2
\scirc u'$ and $u = u' - u'' + \wtilde{y}_0$ where the three summands
have their parameters in $[(i - 1) / 3, i / 3] \times I^{m -1}$, $i =
1, 2, 3$. We lift $u$ to a map $\widehat u = (u, \widehat u_2)$ into
the fiber $F$ of $\proj_2$ (cf \eqref{7.3}) by putting
$$\widehat u_2 (x_1, x') = \proj_2 \scirc u \scirc
  (\text{straight path from } (x_1, x') \text{ to } (0, x') \text{ in } I^m)$$
whenever $(x_1, x') \in \bigl[0, \frac23 \bigr] \times I^{m -1}$ and by using
the strip $\bigl[\frac23, 1\bigr] \times I^{m -1}$ for the obvious deformation
to make sure that $\widehat u$ is constant on $\{ 1\} \times I^{m -1}$
(and hence on the whole boundary $\partial I^m$ of $I^m$). Then $\quot
\scirc \widehat u$ (cf \fullref{7.4}) represents $H ([\wtilde{f}]) := pr_1
\scirc \pinch_* ([\wtilde{f}])$ (cf \eqref{7.5}).

We compare this homotopy class to the degree $\deg^\# (f) = \omega^\#
(f, *)$ of $f = p \scirc \wtilde{f}$ by inspecting the corresponding
geometric (Pontryagin--Thom) data (as described eg in \eqref{2.13} and
\eqref{2.15}). We may assume that $* \in N$ is a regular value of $f$. Applying
the Pontryagin--Thom procedure to $\quot \scirc \widehat u$ we obtain
in the first place the submanifold
$$u^{- 1} ( \{ \wtilde *_1, \dots, \wtilde *_k \} ) \subset
  \mathring{I}{^m} \approx S^m - \{x_0 \}$$
which is the inverse image of the zero section $\Omega (N, y_0)
\times \{ 0\}$ in the Thom space $S^n \wedge (\Omega (N, y_o)^+)$
of the trivial $n$--plane bundle over $\Omega (N, y_0)$ (cf the
discussion preceding \fullref{7.4}). Up to the dilation $\bigl(0,
\frac13\bigr) \times \mathring{I}{^{m -1}} \approx \mathring{I}{^m}$
this submanifold is equal to
$$
u^{'-1} (\{ \wtilde *_1, \dots, \wtilde *_k \} ) = \wtilde
f^{-1} (\{ \wtilde *_1, \dots, \wtilde *_k\} ) = f^{- 1}
(\{ *\} ) = C (f, *) ,
$$
ie to the first component of the triple $(C (f, *), \wtilde{g}_1,
\wbar g^\#)$ which represents $\deg^\# (f)$. Also the maps into
$\Omega (N, y_0)$ are homotopic (up to reversing the direction of
the loops). Indeed, the straight paths which occur in the definition
of $\widehat u_2$ and end in $\{ 0\} \times I^{m -1}$ correspond to a
homotopy $G$ which shifts $C (f, *)$ towards $x_0$ as in the construction
of $\wtilde{g}_1$ (cf \eqref{2.13}); note that the identity map on $\wtilde
N$ can be deformed into the composite
$$
\wtilde{N} \longrightarrow \wtilde{N} / \bigvee^k B^n =
\wtilde{N} - \bigcup^k \mathring{B}{^n} / \bigvee^k S^{n -1}
\cong \wtilde{N}
$$
(the homeomorphisms to the right hand side here and in \eqref{7.1}
agree). However, the framings of the two Pontryagin--Thom data differ by
an automorphism of the trivial $n$--plane bundle induced by $\wtilde
g_1$ (cf \eqref{2.13} and \eqref{2.15}).
\end{proof}

There is a canonical isomorphism
\begin{equation}                                    
\label{7.7}
\kappa \co  \pi_m (S^n \wedge (\Omega N)^+) \longrightarrow \pi_m
\Bigl(\bigvee^k S^n \vee \wtilde{N}, \wtilde{N}\Bigr)
\end{equation}
since the natural epimorphisms from $\pi_m (\vee S^n \vee \wtilde{N})$
onto these groups have the same kernel (cf \eqref{7.5}). The isomorphisms $p_*$
and $\kappa$ allow a purely homotopy theoretical interpretation of our
basic coincidence invariant $\deg^\#$ (cf also \fullref{7.17}).

\begin{cor}         
\label{7.8}
$\kappa \scirc \deg^\# \scirc p_*$ coincides with the composed
homomorphism
$$
\begin{CD}
\pi_m (\wtilde{N}) \longrightarrow \pi_m \Bigl(\wtilde{N}, \wtilde{N} -
\bigcup^k \mathring{B}{^n}\Bigr) @>{ \pinch_* }>> \pi_m \Bigl(\bigvee S^n
\vee \wtilde{N}, \wtilde{N}\Bigr)
\end{CD}
$$
up to an involution of the target group.

In particular, the group $X_m (N)$ (cf \fullref{6.4}) is canonically isomorphic
to the kernel of
$$
(\pinch_*, \partial) \co  \pi_m \Bigl(\wtilde{N}, \wtilde{N} -
\bigcup^k \mathring{B}{^n}\Bigr) \longrightarrow \pi_m \Bigl(\bigvee^k
S^n \vee \wtilde{N}, \wtilde{N}\Bigr) \oplus \pi_{m -1} \Bigl(\wtilde{N} -
\bigcup^k \mathring{B}{^n}\Bigr)
$$
where $\partial$ denotes the obvious connecting homomorphism.
\end{cor}

Next observe that both homomorphisms $i_{1 *}$ and $\pr_1$ in diagram
\eqref{7.5} commute with $\In_*$, $\ret_*$ (cf \eqref{6.15}) on one side and with the
homomorphisms induced by the obvious maps
\begin{equation}                                  
\label{7.9}
  \xymatrix{
    \bigvee S^n \quad\ar@<0.5ex>[r]^-{\operatorname{incl}_1} &
    \quad\bigvee S^n \vee \wtilde{N}
    \ar@<0.5ex>[l]^-{\operatorname{proj}_1}
  }
  \end{equation}
on the other side. Thus we obtain a canonical decomposition
\begin{equation}                                        
\label{7.10}
\pi_m (S^n \wedge (\Omega (N, y_0)^+)) \cong \pi_m \Bigl(\bigvee^k
S^n\Bigr)
\oplus \pi_m \Bigl(\Bigl(\bigvee ^k S^n\Bigr)  \flat \wtilde{N}\Bigr)
\end{equation}
 where $\bigl(\bigvee^k S^n\bigr)  \flat \wtilde{N}$ denotes the homotopy fiber of
 the inclusion
$$\bigvee^k S^n \vee \wtilde{N} \subset \Bigl(\bigvee^k S^n\Bigr)\times \wtilde{N}$$
(cf Ganea \cite[(9)]{G} or
Cornea--Lupton--Oprea--Tanr\'e \cite[6.7]{CLOT}). Given
 $[\wtilde{f}] \in \pi_m (\wtilde{N}, \wtilde{y}_0)$, this yields
 the decomposition
\begin{equation}                                           
\label{7.11}
H ([\wtilde{f}]) = \pr_1 \scirc \pinch_* ([\wtilde{f}]) = ([\coll
\scirc \wtilde{f}], H_{\mathcal C} (\wtilde{f}))
\end{equation}
(cf \eqref{7.5} and \fullref{7.6}) where
$$
\coll \co  \wtilde{N} \longrightarrow \wtilde{N} / \Bigl(\wtilde
N - \bigcup^k \mathring{B}{^n}\Bigr) \cong \bigvee^k S^n
$$
denotes the collapsing map \eqref{7.1} and $H_{\mathcal C} (\wtilde
f)$ is the $\mathcal C$--{\it Hopf--Ganea invariant of} $\wtilde{f}$
based on the cofibration
\begin{equation}                                                
\label{7.12}
\mathcal C \co  \bigvee^k_{i = 1} S^{n -1} \subset
\wtilde{N} - \bigcup^k_{i = 1} \mathring{B}{_i}
\longrightarrow \wtilde{N}
\end{equation}
(cf Cornea--Lupton--Oprea--Tanr\'e \cite[6.7]{CLOT} or
Fern\'andez-Su\'arez--G\'omez-Tato--Tanr\'e \cite[1.1]{FGT}).

\begin{cor}                          
\label{7.13}
Given maps $f, f_1, f_2 \co  (S^m, x_0) \to (N, y_0)$, let $\wtilde
f, \wtilde{f}_1, \wtilde{f}_2$ be the corresponding (basepoint
preserving) liftings to the universal covering space $\wtilde{N}$.

If $\MC (f, *) < \infty$ then $H_{\mathcal C} (\wtilde{f}) = 0$.

If $\MC (f_1, f_2) < \infty$ then $H_{\mathcal C} (\wtilde{f}_1) = -
\underline\inv (H_{\mathcal C} (\wtilde{f}_2))$ where the involution
$\underline\inv$ of
$$
\pi_m \Bigl(\Bigl(\bigvee S^n\Bigr) \flat \wtilde{N}\Bigr) \cong \pi_m (S^n \wedge
(\Omega N)^+) / \In_* \Bigl(\pi_m \Bigl(\bigvee S^n\Bigr)\Bigr)
$$
is induced by the involutions in \eqref{2.12} and \fullref{7.6}.
\end{cor}

\begin{proof} The involution occurring in \fullref{7.6} (and described
at the end of its proof) preserves the subgroup  $\In_* \bigl(\pi_m \bigl(\bigvee
S^n\bigr)\bigr)$ of $\pi_m (S^n \wedge (\Omega N)^+)$ (cf \eqref{6.15}). Thus $H_{\mathcal
C} (\wtilde{f})$ coincides (up to an isomorphism) with the class of
$\deg^\# (f)$ in the quotient group $\pi_m (S^n \wedge (\Omega N)^+) /
\In_* \bigl(\pi_m \bigl(\bigvee S^n\bigr)\bigr)$. But clearly this class vanishes if $\MC (f,
*) < \infty$ (see \fullref{3.5}(i) and \eqref{6.16}). The second claim follows similarly
from \fullref{6.1}.
\end{proof}

\begin{remark}                           
\label{7.14}
Our interpretation of the Hopf--Ganea invariant $H_{\mathcal C}$ as an
obstruction can be extended considerably: given any natural number $k$,
an embedding of $\bigvee^k B^n$ into any simply connected $n$--manifold
$\wtilde{N}$, and any map $\wtilde{f} \co  S^m \to \wtilde{N}$,
$H_{\mathcal C} (\wtilde{f})$  (cf \eqref{7.12}) must vanish if $\wtilde
f$ is homotopic to a map with at least $k$ almost injective points (cf
\fullref{6.22}).
\end{remark}

In view of the last corollary it is natural to ask whether the Hopf--Ganea
invariant is the only finiteness obstruction. In order to get a partial
answer consider the commuting diagram (where $m \ge 2$ as before)
\begin{equation}                                          
\label{7.15}
  \xymatrix{
     \pi_m(N) &
     \pi_m(\wtilde{N})
     \ar[l]^-{\cong}_-{p_*}
     \ar[d]^-{j_*}
     \ar[dr]^-{H_{\mathcal C}}
     \\
     \pi_{m-1}\bigl(\bigvee\limits^{k}_{~} S^{n-1}\bigr)
     \ar[r]_-{E'}
     &
     \pi_m\bigl(\wtilde{N},\wtilde{N} - \bigcup\limits_{~}
     \mathring{B}{^n}\bigr)
     \ar[r]_-{\wtilde{H}_{\mathcal C}}
     &
     \pi_m\bigl(\bigl(\bigvee\limits^{k}_{~} S^n\bigr)\flat\wtilde{N}\bigr)
  }
  \end{equation}
(compare Fern\'andez-Su\'arez--G\'omez-Tato--Tanr\'e
\cite[Section~1]{FGT}). Here $E'$ denotes the \lq
suspension homomorphism\rq\ described by straight paths in $\bigl(\bigvee
B^n, \bigvee S^{n -1}\bigr) \subset \bigl(\smash{\wtilde{N}},
\smash{\wtilde{N}} - \cup \smash{\mathring{B}{^n}}\bigr)$ and $\wtilde{H}_{\mathcal C}$ is a canonical
extension of the $\mathcal C$--Hopf--Ganea homomorphism (cf \eqref{7.12}) to
the indicated relative homotopy group. Given $[\wtilde{f}] \in \pi_m
(\wtilde{N})$, we see (as in the proof of \fullref{6.17}) that $\MC
(p \scirc \wtilde{f}, *) < \infty$ if and only if $j_* ([\wtilde
f])$ lies in $E' (\oplus^k \pi_{m -1} (S^{n -1}))$ (provided $m \ge 3$
or $\pi_1 (N) = 0$).

\begin{thm}                              
\label{7.16}
Assume that $N$ is $q$--connected and $2 \le m \le q + 2n - 3$. Then we
have for all $[\wtilde{f}] \in \pi_m (\wtilde{N})$ :
$$
\MC (p \scirc \wtilde{f}, *) \text{ is finite if and only if }
H_{\mathcal C} (\wtilde{f}) = 0 .
$$
Assume in addition that $E \co  \pi_{m -1} (S^{n -1}) \longrightarrow \pi_m
(S^n)$ is injective. Then $X_m (N) = 0$.
\end{thm}

\begin{proof} When $n = 2 < m$ then $N$ must be open here and $\MC (p
\scirc \wtilde{f}, *) = 0$ (cf \fullref{6.3}). Thus we may assume that $n \ge
3$. Since $\wtilde{N}$ is at least $(m - 2n + 3)$--connected, so is the
(first) inclusion in $\mathcal C$ (cf \eqref{7.12}; use Whitehead's theorem and
excision in homology). Therefore the horizontal line in diagram
\eqref{7.15} is part of an exact EHP--sequence (cf
\cite[1.3]{FGT}). Moreover $E'
(\oplus^k \pi_{m -1} (S^{n -1}))$ is the full image of $E'$ since $k =
\# \pi_1 (N) = 1$ or else $m - 1 < 2n - 3$ (compare Hilton \cite{H}). Our claim
follows from the previous discussion and from \fullref{1.26}.
\end{proof}

\begin{remark}                 
\label{7.17}
We may also consider the $\mathcal C'$--Hopf--Ganea homomorphism
\begin{equation}                            
\label{7.18}
H_{\mathcal C'} \co  \pi_m (N) \longrightarrow \pi_m (S^n \flat N)
\end{equation}
based on the cofibration
\begin{equation}                             
\label{7.19}
\mathcal C'  \co  S^{n -1} \subset N - \mathring{B}{^n}
\longrightarrow  N
\end{equation}
where $S^{n -1}$ is the boundary of some embedded $n$--ball $B^n$
(\lq\lq top cell\rq\rq ) in $N$ (compare \eqref{7.11} and \eqref{7.12}). $H_{\mathcal C'}$
is induced by the corresponding pinch map
\begin{equation}                              
\label{7.20}
\pinch' \co  N \longrightarrow N / S^{n -1} \approx
S^n \vee N
\end{equation}
and the second projection in the canonical decomposition
\begin{equation}                                
\label{7.21}
\pi_m (S^n \vee N) \cong \pi_m (S^n \wedge (\Omega N)^+) \oplus \pi_m
(N) \cong \pi_m (S^n) \oplus \pi_m (S^n \flat N) \oplus \pi_m (N)
\end{equation}
(compare \eqref{7.1}, \eqref{7.5}, and \eqref{7.10}); here the first component homomorphism
of $\pinch'_*$ is induced by a collapsing map $\coll' \co  N \to S^n$
of degree $\pm 1$ (as in \eqref{2.16}).

According to \fullref{7.6} the combined homomorphism
\begin{equation}                                
\label{7.22}
(\coll'_*,  H_{\mathcal C'}) \co  \pi_m (N) \longrightarrow \pi_m
(S^n) \oplus \pi_m (S^n \flat N)
\end{equation}
agrees with $\deg^\#$ up to an isomorphism between the target groups. In
particular, for all $[f] \in \pi_m (N)$ the Hopf--Ganea invariant
$H_{\mathcal C'} ([f])$ (relative to the attaching class of the top
cell of $N$) can be described entirely by the coincidence data of the
preimage $f^{-1} (\{ *\})$. To be more precise assume that $f$ is smooth
with regular value $* \in N$. Then $\coll'_* ([f])$ corresponds (via
Pontryagin--Thom) to the bordism class of the framed submanifold $C :=
C (f, *) = f^{- 1} (\{ *\} )$ of $S^m - \{ x_0\}$ and $H_{\mathcal C'}
([f])$ measures the added information which the lifting $\wtilde{g}$
(cf \eqref{1.6}) or, equivalently, the map $\wtilde{g}_1 \co  C \to \Omega
N$ contributes to $\deg^\# ([f]) = [f^{- 1} (\{ *\} ), \wtilde
g_1, \wbar g^\#_1 ]$ (cf \eqref{2.13}--\eqref{2.15}). On the other hand,
the element $\coll_* ([\wtilde{f}]) \in \pi_m (\vee^k S^m)$ (cf
\eqref{7.10}, \eqref{7.11}) corresponds to the link $C = \amalg C_\alpha$
of disjoint framed submanifolds (ie it captures the Nielsen decomposition of $C$
determined by $\wtilde{g}$, cf \eqref{2.6}), and the Ganea--Hopf invariant
$H_{\mathcal C} ([\wtilde{f}])$, which we studied in
\eqref{7.11}--\eqref{7.15}, \fullref{7.13} and \fullref{7.16},
measures the remaining information contained in $\deg^\# ([f])$.

In general, $H_{\mathcal C'} ([f])$ contains much more information than
$H_{\mathcal C} ([\wtilde{f}])$. Indeed,
\begin{equation}                                 
\label{7.23}
\pi_m (S^n \flat N) \cong \ker (\ret'_*) \oplus \pi_m
\Bigl(\Bigl(\bigvee^k S^n\Bigr)
\flat \wtilde{N}\Bigr)
\end{equation}
where $\ret' \co  \vee^k S^n \to S^n$ maps each sphere in the wedge
identically to $S^n$. In view of the Hilton decomposition
\begin{equation}                                 
\label{7.24}
\pi_m \Bigl(\bigvee^k S^n\Bigr) \cong \underset{k}{\bigoplus} \pi_m (S^n)
\oplus \underset{k \choose 2}{\bigoplus}  \pi_m (S^{2n -1}) \oplus \dots
\end{equation}
(cf Hilton \cite{H}) the kernel of $\ret'_*$ in \eqref{7.23} may be highly nontrivial
whenever $k > 1$.

In particular, $H_{\mathcal C'} ([f])$ need not vanish when $\MC (f, *)$
is finite (eg in the case $m = n$ and $\# G > 1$ in \fullref{1.25}). Thus
it is suitable to use the universal covering space $\wtilde{N}$ and the
corresponding Hopf--Ganea homomorphism $H_{\mathcal C} ([\wtilde{f}])$
when we develop such finiteness criteria as \fullref{7.13} and
\fullref{7.16}. However,
$\wtilde{N}$ is not required in \fullref{7.8}: clearly $X_m (N)$ is also
canonically isomorphic to the kernel of the obvious homomorphism
$$
(\pinch'_*, \partial') \co  \pi_m (N, N - \mathring{B}{^n})
\longrightarrow \pi_m (S^n \vee N, N) \oplus \pi_{m -1} (N -
\mathring{B}{^n}) .
$$
\end{remark}

\appendix

\section{Base points}                    
\setobjecttype{App}
\label{secA}

It is sometimes useful to require that the maps $f_1$ and $f_2$ as well
as their homotopies preserve base points. As we will see this has no
impact on our numerical invariants whenever $m, n \ge 1$.

Let $x_0 \in M$ and $y_1 \ne y_2 \in N$ be a given choice of base
points. For any two \lq\lq base point preserving\rq\rq\ maps
$$f_i \co  (M, x_0) \longrightarrow (N, y_i),\quad i = 1, 2,$$
define $\MCC_b (f_1, f_2)$ (and $\MC_b (f_1, f_2)$, resp.) to be
the minimum number of path components (and of points, resp.) in the
coincidence locus $C (f'_1, f'_2)$ of any pair of maps in the same base
point preserving homotopy classes, ie
$$[f'_i] = [f_i] \in [(M, x_0), (N, y_i)],\quad  i = 1, 2 .$$

\begin{lemA}  $\MCC_b (f_1, f_2) = \MCC (f_1, f_2)$ and $\MC_b (f_1, f_2)
= \MC (f_1, f_2)$; in particular if the pair $(f_1, f_2)$ is loose then
there are already base point preserving homotopies which deform $f_1$
and $f_2$ away from one another.
\end{lemA}

(For a related result in the context of classical fixed point theory
see eg Jiang \cite[Section~3]{Ji1}).

\begin{proof} The claims concerning $\MCC$ and $\MC$ will be treated
simultaneously. We may concentrate on  the case $m, n \ge 2$.

Clearly $M(C)C_b (f_1, f_2) \ge M(C)C (f_1, f_2)$.

Conversely, let the minimum number $M(C)C (f_1, f_2)$ be realized by
a pair $(f''_1, f''_2)$. Without changing the number of coincidence
components (or points, resp.) we may -- in a first step -- deform this
pair until it preserves base points. Indeed, if there is a point $\widehat
x_0 \in M$ such that $f''_1 (\widehat x_0) \ne f''_2 (\widehat x_0)$,
compose $(f''_1, f''_2)$ with isotopies in $M$ and $N$ which move $x_0$
to $\widehat x_0$ and $f_i (\widehat x_0)$ to $y_i$, $i = 1, 2$; if $f''_1
\equiv f''_2$ remove first a small ball $\mathring{B} \subset M$
from the coincidence set $C (f''_1, f''_2)$ by \lq\lq pushing $f''_2$
slightly off $f''_1$ in $\mathring{B}$\rq\rq.

The second step is a simple modification near $x_0$ which makes $f''_i$
homotopic to $f_i$ in the base point preserving sense, $i = 1, 2$.  Let
$$
c_i = H_i (x_0, -) \co  (I, \{0, 1\}) \longrightarrow (N, y_i)
$$
be smooth transverse paths resulting from homotopies $H_i \co  f_i \sim
f''_i$, $i = 1, 2$. After suitable changes of the parametrization near
finitely many intersection points $c_1$ and $c_2$ will be coincidence
free. Now identify a small ball in $M$ around $x_0$ with the unit ball
$B^m$ in $\mathbb R^m$. We may assume that for every $x \in B^m$, $f''_i
(x) = y_i$; then replace this by $f''_i (x) := c_i (\Vert x \Vert)$, $i
= 1, 2$. This procedure does not change the coincidence set of $(f''_1,
f''_2)$ but yields basepoint preserving homotopies. Thus in the end we
conclude that $M(C)C_b (f_1, f_2) \le \# (\pi_0) C (f''_1, f''_2)
= M(C)C (f_1, f_2)$.  \end{proof}

\bibliographystyle{gtart}
\bibliography{link}

\end{document}